\documentclass[11pt, one side,emlines]{amsart}
\usepackage{tikz}
\usepackage{amssymb,latexsym,xy,eucal,mathrsfs}
\usepackage{bm,amsfonts,amssymb,amsthm,newlfont,eucal, mathrsfs,latexsym, caption}
\usetikzlibrary{decorations.pathmorphing}

\usepackage{float}
\usepackage{ circuitikz }
\usepackage{algorithm}
\usepackage{algorithmic}
\textwidth=15.2cm \textheight=24cm \theoremstyle{plain}

\theoremstyle{definition}

\newtheorem{theorem}{Theorem}

\newtheorem{lemma}{Lemma}

\newtheorem{corollary}{Corollary}

\numberwithin{equation}{section}  \voffset
-55truept \hoffset -20truept
\begin{document}
	\baselineskip 15truept
	\title[The 3-path-connectivity of the augmented cubes]{The 3-path-connectivity of the augmented cubes}
	
	%	\subjclass[2010]{Primary 16W10; Secondary 06A06; 47L30} 
	\maketitle 
		\begin{center} 
	\author{S. A. Kandekar, R. Barabde, S. A. Mane }\\	
					{\small Center for Advanced Studies in Mathematics,
						Department of Mathematics,\\ Savitribai Phule Pune University, Pune-411007, India.}\\
				\email{\emph{smitakandekar54@gmail.com; rushikeshbarbde@gmail.com; manesmruti@yahoo.com}} 
			\end{center}
	\begin{abstract}

		Connectivity is a cornerstone concept in graph theory, essential for evaluating the robustness of networks against failures. To better capture fault tolerance in complex systems, researchers have extended classical connectivity notions—one such extension being the $k$-path-connectivity, $\pi_k(G)$, introduced by Hager. Given a connected simple graph $G = (V, E)$ and a subset $D \subseteq V$ with $|D| \geq 2$, a $D$-path is a path that includes all vertices in $D$. A collection of such paths is internally disjoint if they intersect only at the vertices of $D$ and share no edges. The maximum number of internally disjoint $D$-paths in $G$ is denoted $\pi_G(D)$, and the $k$-path-connectivity is defined as
		
		$$
		\pi_k(G) = \min \{ \pi_G(D) \mid D \subseteq V(G),\ |D| = k \}.
		$$
		
		In this paper, we investigate the 3-path-connectivity of the augmented cube $AQ_n$, a variant of the hypercube known for its enhanced symmetry and fault-tolerant structure. We establish the exact value of $\pi_3(AQ_n)$ and show that:
		
		$$
		\pi_3(AQ_n) = 
		\begin{cases}
			\frac{3n}{2} - 2, & \text{if } n \text{ is even}, \\
			\frac{3(n - 1)}{2} - 1, & \text{if } n \text{ is odd}.
		\end{cases}
		$$

	\end{abstract}
	\section{Introduction}
	
	Let $G = (V, E)$ be a connected simple graph, and let $D \subseteq V$ be a subset with at least two vertices. A path $P$ in $G$ is called a $D-$path if all the vertices in $D$ lie on $P$. A collection of paths $P_1, P_2, \dots, P_m$ is called a set of internally disjoint $D-$paths if, for all $1 \leq i \neq j \leq m$, the only common vertices between $P_i$ and $P_j$ are those in $D$, and they share no edges; formally,
	$V(P_i) \cap V(P_j) = D \quad \text{and} \quad E(P_i) \cap E(P_j) = \emptyset.$
	
	Let $\pi_G(D)$ be the maximum number of internally disjoint $D-$paths in $G$. For any integer $k$ with $2 \leq k \leq |V(G)|$, the \textit{k-path-connectivity} of $G$, written $\pi_k(G)$, is defined as: $\pi_k(G) = \min\{ \pi_G(D) \mid D \subseteq V(G), |D| = k \}.$
	
	Graph connectivity is a fundamental concept used to assess network's resilience against failures. To more accurately reflect the fault-tolerant capabilities of networks, several generalizations of classical connectivity have been proposed. One such generalization is the $k$-path-connectivity, denoted by $\pi_k(G)$, introduced by Hager \cite{hager}. Another important extension is the generalized connectivity, $\kappa_k(G)$, introduced by Chartrand et al. \cite{chartrand1984generalized}, which broadens the classical ideas of connectivity and path-connectivity. This parameter has received considerable attention in the literature, with numerous related works exploring its properties and applications. For more detailed discussions, see \cite{generalizedhypercube}, \cite{generalizedfolded}, \cite{mane2024pendant}, \cite{cheng2023generalized}.
	
	Li et al. \cite{li2021tree} established that, for any fixed integer $k \geq 1$, determining whether $\pi_G(D) \geq k$ for a given subset $D \subseteq V(G)$ is an NP-complete problem, highlighting the computational complexity of $k$-path-connectivity.
	
	In particular, the 3-path-connectivity $\pi_3(G)$ has been extensively investigated, and exact values have been determined for several well-known graph classes. These include the complete graph $K_n$ and the complete bipartite graph $K_{s,t}$ \cite{hager}, the hypercube $Q_n$ \cite{zhu}, the $k$-ary $n$-cube $Q_n^k$ \cite{C2}, the pancake graph \cite{wang}, and the folded hypercube \cite{wang20253} . Furthermore, tight bounds for $\pi_3(G)$ have been derived in cases where the graph $G$ is constructed via the lexicographic product of two graphs \cite{mao}, demonstrating the richness and depth of this parameter across various graph operations and structures.
	
	The augmented cube $AQ_n$ is a variation of the classical hypercube that possesses several enhanced structural properties, making it a valuable topology for interconnection networks. 
	In this paper, by closely analyzing the structure of $(AQ_n)$ we investigate the 3-path connectivity of augmented cubes $\pi_3(AQ_n)$. We prove that 3-path connectivity of the augmented cube is\\
	      $$
	      \pi_3(AQ_n) = 
	      \begin{cases}
	      	\frac{3n}{2} - 2, & \text{if } n \text{ is even}, \\
	      	\frac{3(n - 1)}{2} - 1, & \text{if } n \text{ is odd}.
	      \end{cases}
	      $$

	This result contributes to the understanding of fault-tolerant properties and path diversity in augmented cube networks.
	For additional terminology and notation, refer to \cite{west2001}.

	\section{\textbf{Preliminaries}}
	
	In this section, we begin by recalling the definition of the augmented cube, originally introduced by Choudum et al. \cite{choudum}. We then introduce the notation used in the decomposition of the augmented cube. While the augmented cube is commonly divided into two parts, we adopt a four-part decomposition. This refined approach not only clarifies the structure but also results in a shorter and more straightforward proof.
	
	The augmented cube of dimension \(n\), denoted \(AQ_n\), is a graph with vertex set 
	\(V(AQ_n) = \{x_1x_2 \dots x_n \mid x_i \in \{0, 1\}, 1 \leq i \leq n\}\), consisting of \(2^n\) vertices. 
	For \(n = 2\), the augmented cube is \(AQ_2 = K_2\), a complete graph with two vertices. 
	For \(n \geq 2\), \(AQ_n\) is constructed recursively from two disjoint copies of \(AQ_{n-1}\), denoted 
	\(AQ_{n-1}^0\) and \(AQ_{n-1}^1\), with the vertices of \(AQ_{n-1}^0\) and \(AQ_{n-1}^1\) prefixed by \(0\) and \(1\), respectively. 
	The vertex set of \(AQ_n\) is therefore \(V(AQ_n) = V(AQ_{n-1}^0) \cup V(AQ_{n-1}^1)\), and its edge set includes all edges within \(AQ_{n-1}^0\) and \(AQ_{n-1}^1\), along with hypercubic edges connecting each vertex \(x = 0x_1x_2 \dots x_{n-1} \in V(AQ_{n-1}^0)\) to \(x^h = 1x_1x_2 \dots x_{n-1} \in V(AQ_{n-1}^1)\) and complementary edges connecting \(x\) to \(x^c = 1\overline{x_1x_2 \dots x_{n-1}} \in V(AQ_{n-1}^1)\), where \(\overline{x_1x_2 \dots x_{n-1}}\) denotes the bitwise complement of \(x_1x_2 \dots x_{n-1}\). 
	
	The augmented cube \(AQ_{n-1}^i\) of dimension \(n-1\) is similarly formed by two copies of \(AQ_{n-2}\), denoted \(AQ_{n-2}^{ij}\) and \(AQ_{n-2}^{i(1-j)}\) $(0 \leq i, j \leq 1)$. 
	The augmented cube \(AQ_n\) can be decomposed into four copies of augmented cubes of dimension \(n-2\) as 
	$$AQ_n = (AQ_{n-2}^{00} \square AQ_{n-2}^{01}) \square (AQ_{n-2}^{10} \square AQ_{n-2}^{11})$$ Where $AQ_{n-1}^0\square AQ_{n-1}^1$ is recursive construction of $AQ_n$ as described above.
	
	There are two perfect matchings between the subgraphs $AQ_{n-2}^{ij}$ and $AQ_{n-2}^{i(1-j)}$:
	
	1. \textit{Hypercubic matching}: Each vertex $x \in AQ_{n-2}^{ij}$ is connected to its \textit{hypercubic neighbor} $x^{h_{n-1}} \in AQ_{n-2}^{i(1-j)}$.
	
	2. \textit{Complementary matching}: Each vertex $x \in AQ_{n-2}^{ij}$ is connected to its \textit{complementary neighbor} $x^{c_{n-1}} \in AQ_{n-2}^{i(1-j)}$.
	
	These matchings occur within the subgraph $AQ_{n-1}^i = AQ_{n-2}^{ij} \square AQ_{n-2}^{i(1-j)}$.
	
	Similarly, in the full graph $AQ_n = AQ_{n-1}^0 \square AQ_{n-1}^1$, we have:
	
	3. A perfect \textit{hypercubic matching} between $AQ_{n-2}^{ij}$ and $AQ_{n-2}^{(1-i)j}$, where each vertex $x \in AQ_{n-2}^{ij}$ is connected to its hypercubic neighbor $x^h \in AQ_{n-2}^{(1-i)j}$.
	
	4. A perfect \textit{Complementary matching} between $AQ_{n-2}^{ij}$ and $AQ_{n-2}^{(1-i)(1-j)}$, where each vertex $x \in AQ_{n-2}^{ij}$ is connected to its complementary neighbor $x^c \in AQ_{n-2}^{(1-i)(1-j)}$.\\
	
	Let $AQ_{n-2}^{ij}\diamond AQ_{n-2}^{(1-i)j}$ be a subgraph of $AQ_n$ containing $AQ_{n-2}^{ij},AQ_{n-2}^{(1-i)j}$ and the perfect matching between them for $0\leq i,j\leq1$. Also, let $P(x,y)$ be a path in between the vertex $x$ and $y$ in a graph.
	
	\section*{\textbf{The 3-Path connectivity of augmented cubes}}

	In this section, we begin by recalling several lemmas established in previous works, which will serve as essential tools in the proof of our main theorem. We then present a proof of the 3-path connectivity of the augmented cube of dimension $4$, and subsequently extend the result to augmented cubes of dimension $n$ using mathematical induction. While the augmented cube is traditionally decomposed into two parts, we employ a four-part decomposition. This refined approach not only offers a clearer understanding of the structure but also leads to a more concise and elegant proof.
	
	\begin{lemma}\cite{west2001}\label{lemma1}
		If $G$ is a simple $k$-connected graph, then for any vertex $x\in V(G)$ and any subset $S\subseteq V(G)$ with $|S|=k$, there exists $k-$vertex disjoint paths in $G$, each joining $x$ to every vertex in $S$.
	\end{lemma}
	\begin{lemma}\cite{west2001}\label{lemma7}
		Let $G$ be a simple $k$-connected graph, and let $A, B$ be subsets of $V(G)$ with $|A| = |B| = k$. Then there exist $k$ vertex-disjoint paths in $G$, each join a distinct vertex in $A$ to a distinct vertex in $B$.

	\end{lemma}
	\begin{lemma}\cite{west2001}\label{lemma2}
		Let $G_1$ and $G_2$ are two $k-$connected graphs. Suppose there exist a perfect matching $M_1$ that pairs each vertex in $G_1$ with a distinct vertex in $G_2$. Define a graph $G$ as the union of $G_1$ and $G_2$ along with the edges in $M_1$ such that $G=V(G_1)\cup V(G_2)\cup E(G_1)\cup E(G)_2\cup M_1$. Then $G$ is $(k+1)$-connected.
	\end{lemma}
	\begin{lemma}\cite{choudum}\label{lemma3}
		Augmented cube $AQ_n$ of dimension $n$ is $2n-1$ connected, for $n\neq 3$ and $4$ connected for $n=3$. 
	\end{lemma}
	\begin{lemma}\cite{C2}\label{lemma4}
		Let $G$ be a $ k$-regular graph. Then $\pi_3 (G) \leq \lfloor \frac{3k-r}{4}\rfloor$
		, where\\
		$r =\text{max }\{|N_G(x) \cap N_G (y) \cap N_G (z)| : \{x, y, z\} \subseteq V (G)\}$.
	\end{lemma}
	
	\begin{lemma}\cite{C1} \label{lemma5} Any two distinct vertices in $AQ_n$ for $n \geq 3$, have at most four common neighbors.
	\end{lemma}
	
	\begin{lemma} \label{lemma6}
		
		For the augmented cube $AQ_n$ of dimension $n \geq 4$, the 3-path-connectivity satisfies:
		
		$$
		\pi_3(AQ_n) \leq
		\begin{cases}
			\frac{3n}{2} - 2, & \text{if } n \text{ is even}, \\
			\frac{3(n - 1)}{2} - 1, & \text{if } n \text{ is odd}.
		\end{cases}
		$$

	\end{lemma} 
	\begin{proof}	By Lemma \ref{lemma5}, any two distinct vertices in the augmented cube $AQ_n$ have at most four common neighbors. Consider the vertices
		
		$$
		x = 000x_1x_2\dots x_{n-3}, \quad y = 011\overline{x_1x_2\dots x_{n-3}}, \quad z = 101x_1x_2\dots x_{n-3},
		$$
		
		where $x_1, x_2, \dots, x_{n-3} \in \{0,1\}$. Then the vertices $x, y,$ and $z$ have exactly four common neighbors, namely:
		
		$$
		\begin{aligned}
			a &= 001\overline{x_1x_2\dots x_{n-3}}, \\
			b &= 111\overline{x_1x_2\dots x_{n-3}}, \\
			c &= 100x_1x_2\dots x_{n-3}, \\
			d &= 010x_1x_2\dots x_{n-3},
		\end{aligned}
		$$
		
		which is the maximum possible by Lemma \ref{lemma5}. Therefore, for $AQ_n$, we have $r = 4$, and it follows that
		
		$$
		\pi_3(AQ_n) \leq \left\lfloor \frac{3(2n - 1) - 4}{4} \right\rfloor = \left\lfloor \frac{6n - 3}{4} \right\rfloor - 1.
		$$
		
		In particular,
		
		 if $n$ is even, then $\pi_3(AQ_n) \leq \frac{3n}{2} - 2$;\\
		and  if $n$ is odd, then $\pi_3(AQ_n) \leq \frac{3(n - 1)}{2} - 1$.

	\end{proof}
	\begin{lemma} \label{lemma8}
		
		For the augmented cube $AQ_4$, we have 
			$\pi _3(AQ_4)=4$
	\end{lemma}
	\begin{proof}
		\textbf{Case 1 :} $\{x,y,z\}\subseteq V(AQ_{2}^{ij})$ for some $i$ and $j$, $0\leq i,j\leq 1$.\\
		
		Without loss of generality consider $\{x,y,z\}\subset AQ_2^{00}$.
		Let, $x=00x_1x_2,\\ y=00\overline{x_1}x_2, z=00x_1\overline{x_2}.$ Then we get four vertex disjoint paths containing $x,y$ and $z$ in $AQ_4$ as follows, (See Figure 1).\\
		
		$\psi_1=\{yx, xz\}\\$
		
		$\psi_2=\{xa, ay, yz\}\\$
		
		$\psi_3=\{xx^{h_{3}}, x^{h_{3}}y^{h_{3}}, y^{h_{3}}y, yz^{h_3}, z^{h_3}z\}\\$
		
		$\psi_4=\{yy^h ,y^hx^h,x^hx,  xx^c, x^cz^c, z^cz
		$\}\\
		\begin{figure}[ht]
			\centering
			\begin{tikzpicture}
				% Define the vertices with smaller size and labels outside
				\draw[thick] (-.7, 3.5) rectangle (2.2, 6);
				\draw[thick] (3.2, 3.5) rectangle (6.1, 6);
				\draw[thick] (-.7, 2.5) rectangle (2.2, 0);
				\draw[thick] (3.2, 2.5) rectangle (6.1, 0);
				\node[circle, draw, fill=blue!20, minimum size=4pt, inner sep=0pt, label=left:z] (z) at (0, 4) {};
				\node[circle, draw, fill=blue!20, minimum size=4pt, inner sep=0pt, label=right:a] (a) at (1.5, 4) {};
				\node[circle, draw, fill=blue!20, minimum size=4pt, inner sep=0pt, label=right:y] (y) at (1.5, 5.5) {};
				\node[circle, draw, fill=blue!20, minimum size=4pt, inner sep=0pt, label=left:x] (x) at (0, 5.5) {};
				\node[circle, draw, fill=blue!20, minimum size=4pt, inner sep=0pt, label=left:$x^{h_3}$] (E) at (4, 5.5) {};
				\node[circle, draw, fill=blue!20, minimum size=4pt, inner sep=0pt, label=right:$y^{h_3}$] (F) at (5.5, 5.5) {};
				\node[circle, draw, fill=blue!20, minimum size=4pt, inner sep=0pt, label=left:$z^{h_3}$] (G) at (4, 4) {};
				\node[circle, draw, fill=blue!20, minimum size=4pt, inner sep=0pt] (H) at (5.5, 4) {};
				\node[circle, draw, fill=blue!20, minimum size=4pt, inner sep=0pt, label=left:$x^h$] (xh) at (0, 2) {};
				\node[circle, draw, fill=blue!20, minimum size=4pt, inner sep=0pt, label=right:$y^h$] (yh) at (1.5, 2) {};
				\node[circle, draw, fill=blue!20, minimum size=4pt, inner sep=0pt] (c) at (4, 2) {};
				\node[circle, draw, fill=blue!20, minimum size=4pt, inner sep=0pt, label=right:$z^c$] (zc) at (5.5, 2) {};
				\node[circle, draw, fill=blue!20, minimum size=4pt, inner sep=0pt] (zh) at (0, 0.5) {};
				\node[circle, draw, fill=blue!20, minimum size=4pt, inner sep=0pt] (f) at (1.5, 0.5) {};
				\node[circle, draw, fill=blue!20, minimum size=4pt, inner sep=0pt] (yc) at (4, 0.5) {};
				\node[circle, draw, fill=blue!20, minimum size=4pt, inner sep=0pt, label=right:$x^c$] (xc) at (5.5, 0.5) {};
				\node at (0.8, 7){\textbf{AQ$_2^{00}$}};
				\node at (4.7, 7){\textbf{AQ$_2^{01}$}};
				\node at (4.7, -0.5){\textbf{AQ$_2^{11}$}};
				\node at (0.8, -0.5){\textbf{AQ$_2^{10}$}};
				% Draw the edges
				\draw [thick, draw=green] (a) to (y);
				\draw [thick, draw=red] (y) -- (x);
				\draw[thick, draw=red] (x) -- (z);
				\draw[thick, draw=green] (x) -- (a);
				\draw[thick, draw=green] (z) -- (y);
				\draw[thick, bend right=50, draw=blue] (z) to (G);
				\draw[thick, draw=blue] (G) -- (y);
				\draw[thick, bend left=30, draw=blue] (x) to (E);
				\draw[thick, bend left=30, draw=blue] (y) to (F);
				\draw[thick, draw=blue] (F) -- (E);
				\draw[thick, bend right=30] (x) to (xh);
				\draw[thick, bend left=30] (y) to (yh);
				\draw[thick] (yh) to (xh);
				\draw[thick, bend right=20] (x) to (xc);
				\draw[thick] (z) to (zc);
				\draw[thick] (xc) to (zc);
				\draw[dotted, line width=.5mm] (-2, 3) -- (7.4, 3);
			\end{tikzpicture}
			\caption{Case 1}
			\label{fig:graph-example}
		\end{figure}

		\noindent\textbf{Case 2 :} $x,y\in V(AQ_{2}^{00})$ and $z\in V(AQ_{2}^{01})$\\
		
		Consider, $x,y\in V(AQ_{2}^{00}\diamond AQ_{2}^{10}).$ By Lemma \ref{lemma1} and Lemma \ref{lemma2}, $AQ_{2}^{00}\diamond AQ_{2}^{10}$ is $4$-connected which gives four vertex disjoint paths joining $x$ and $y$ in $AQ_{2}^{00}\diamond AQ_{2}^{10}$, say, $P_1,P_2,P_3,P_4$. Now, $X=\{x^{h_{3}},y^{h_{3}},x^c,y^c\}$ is the set of adjacent vertices of $x$ and $y$ in $AQ_{2}^{01}\diamond AQ_{2}^{11}$. Since $AQ_2^{01}\diamond AQ_2^{11}$ is four connected, there are four vertex disjoint paths, each joining $z$ to every vertex of set $X$. let, $P(x^{h_{3}},z),P(y^{h_{3}},z),P(x^c,z)$ and $P(y^c,z)$ be the paths between $z$ and every vertex of set $X$ in $AQ_{2}^{01}\diamond AQ_{2}^{11}$. \\
		We get four vertex disjoint paths containing $x, y,$ and $z$ in $AQ_n$ as follows, (See Figure 2).\\
		$\psi_1= P_1\cup \{xx^{h_{3}}\}\cup P(x^{h_{3}},z)\\
		\psi_2=P_2\cup \{yy^{h_{3}}\}\cup P(y^{h_{3}},z)\\
		\psi_3=P_3\cup \{xx^c\}\cup P(x^c,z)\\
		\psi_4=P_4\cup \{yy^c\}\cup P(y^c,z).$\\
		\begin{figure}[H]
			\centering
			\begin{tikzpicture}
				% Draw the first rectangle
				\draw[thick] (0, 0) rectangle (3, 5.5);
				\draw[thick] (-4, 0) rectangle (-1., 5.5);
				
				\draw[dotted] (0.2, 0.5) rectangle (2.8, 2.3);
				\draw[dotted] (0.2, 5.2) rectangle (2.8, 3.3);
				\draw[dotted] (-3.8, 0.5) rectangle (-1.2, 2.3);
				\draw[dotted] (-3.8, 5.2) rectangle (-1.2, 3.3);
				
				% Labels for the rectangles
				\node at (-2.3, -1) {\textbf{$(AQ_2^{00}\diamond AQ_2^{10})$}};
				\node at (1.8, -1) {\textbf{$(AQ_2^{01}\diamond AQ_2^{11})$}};
				
				% Smaller vertices and connections
				\node[circle, draw, fill=blue!20, minimum size=4pt, inner sep=0pt, label=left:$x$] (x) at (-3.3, 5.) {};
				\node[circle, draw, fill=blue!20, minimum size=4pt, inner sep=0pt] (a) at (-3.3, 2) {};
				\node[circle, draw, fill=blue!20, minimum size=4pt, inner sep=0pt] (b) at (-1.6, 0.8) {};
				\node[circle, draw, fill=blue!20, minimum size=4pt, inner sep=0pt, label=right:$y$] (y) at (-1.5, 3.6) {};
				
				% Edges in the left rectangle
				\draw[decorate, draw=red, decoration={snake, amplitude=0.02cm, segment length=0.1cm}, bend left=45] (x) to (y);
				\draw[decorate, draw=blue, decoration={snake, amplitude=0.02cm, segment length=0.1cm}] (x) to (y);
				\draw[decorate, draw=green, decoration={snake, amplitude=0.02cm, segment length=0.1cm}, bend right=45] (x) to (y);
				\draw[decorate, decoration={snake, amplitude=0.02cm, segment length=0.1cm}] (x) to (a);
				\draw[decorate, decoration={snake, amplitude=0.02cm, segment length=0.1cm}] (b) to (y);
				\draw[decorate, decoration={snake, amplitude=0.02cm, segment length=0.1cm}] (a) to (b);
				
				% Smaller vertices in the right rectangle
				\node[circle, draw, fill=blue!20, minimum size=4pt, inner sep=0pt, label=left:$x^{h_3}$] (xh2) at (1, 5) {};
				\node[circle, draw, fill=blue!20, minimum size=4pt, inner sep=0pt, label=left:$y^{h_3}$] (yh2) at (2.5, 3.6) {};
				\node[circle, draw, fill=blue!20, minimum size=4pt, inner sep=0pt, label=left:$y^c$] (yc) at (1, 2) {};
				\node[circle, draw, fill=blue!20, minimum size=4pt, inner sep=0pt, label=left:$x^c$] (xc) at (2.3, 0.8) {};
				\node[circle, draw, fill=blue!20, minimum size=4pt, inner sep=0pt, label=left:$z$] (z) at (1.8, 4.3) {};
				
				% Edges in the right rectangle
				\draw[decorate, draw=red, decoration={snake, amplitude=0.02cm, segment length=0.1cm}] (z) to (xh2);
				\draw[decorate, draw=blue, decoration={snake, amplitude=0.02cm, segment length=0.1cm}] (z) to (yh2);
				\draw[decorate, decoration={snake, amplitude=0.02cm, segment length=0.1cm}] (z) to (yc);
				\draw[decorate, draw=green, decoration={snake, amplitude=0.02cm, segment length=0.1cm}] (z) to (xc);
				\draw[thin,draw=red,bend left=15] (x) to (xh2);
				\draw[thin,draw=blue,bend right=15] (y) to (yh2);
				\draw[thin,draw=green,bend right=15] (x) to (xc);
				\draw[thin] (y) to (yc);
				% Additional labels
				\node at (-1.8, 4.8) {\textbf{$P_1$}};
				\node at (-2.1, 4.3) {\textbf{$P_2$}};
				\node at (-2.5, 3.7) {\textbf{$P_3$}};
				\node at (-3.6, 3.3) {\textbf{$P_4$}};
				
				% Dotted line
				\draw[dotted, line width=.5mm] (-4, 3) -- (3.4, 3);
				
				% Labels for sub-rectangles
				\node at (-2.3, 6.3) {\textbf{$AQ_2^{00}$}};
				\node at (1.9, 6.3) {\textbf{$AQ_2^{01}$}};
				\node at (1.9, -0.3) {\textbf{$AQ_2^{11}$}};
				\node at (-2.3, -0.3) {\textbf{$AQ_2^{10}$}};
			\end{tikzpicture}
			
			\caption{Case 2}
			\label{fig:graph-example}
		\end{figure}

		\textbf{Case 3 :}
		$x,y\in V(AQ_{3}^0)$ and $z\in V(AQ_{3}^1).$\\
		\textbf{Subcase 3.1 :} Suppose $x^{c_3}=y$.\\
		Without loss of generality suppose $x\in V(AQ_2^{00})$ which gives $x^{c_3}=y\in V(AQ_2^{01})$. Let $x_1,x_2,x_3$ be the adjacent vertices of $x$ in $V(AQ_{2}^{00})$ and $y_1,y_2,y_3$ be the adjacent vertices of $y$ in $V(AQ_{2}^{01})$. Without loss of generality assume that $x^{h_{3}}=y_3$ and $y^{h_3}=x_3.$ Also $x_1^{h_3}=y_1$ and $x_2^{h_3}=y_2.$ Then there are five vertex disjoint paths joining $x$ and $y$ in $AQ_{3}^0$ as follows,\\
		$P_1=\{xy\}\\
		P_2=\{xx_3, x_3y\}\\
		P_3=\{xy_3, y_3y\}\\
		P_4=\{xx_1, x_1y_1, y_1y\}\\
		P_5=\{xx_2, x_2y_2, y_2y\}\\$
		Consider $X=\{x^h,y^h,x^h_2,y^h_2\}\subset V(AQ_{3}^1)$. Then there are four vertex disjoint paths each joining $z$ to every vertex in $X,$ say, $P(x^h,z),P(y^h,z),P(x^h_2,z)$ and $(y^h_2,z)$. 
		Thus we get four vertex disjoint paths containing $x,y,$ and $z$ in $AQ_4$ as follows, (See Figure 3).\\
		\begin{figure}[H]
			\centering	
			\begin{tikzpicture}[baseline=(current bounding box.south), every node/.style={inner sep=0, outer sep=0}]
				
				\draw[thick] (-0.5, 6) rectangle (2.7, 3.5);
				\draw[thick] (3.2, 6) rectangle (6.5, 3.5);
				\draw[dotted, line width=.5mm] (-2, 3) -- (7.4, 3);	
				% Define the vertices with labels outside
				\node[circle, draw, fill=blue!20,minimum size=4pt, inner sep=0pt, label=below:$x_2$] (x2) at (0, 4) {};
				\node[circle, draw, fill=blue!20,minimum size=4pt, inner sep=0pt, label=below:$x_3$] (x3) at (1.5, 4) {};
				\node[circle, draw, fill=blue!20,minimum size=4pt, inner sep=0pt, label=above:$x_1$] (x1) at (1.5, 5.5) {};
				\node[circle, draw, fill=blue!20,minimum size=4pt, inner sep=0pt, label=above:$x$] (x) at (0, 5.5) {};
				\node[circle, draw, fill=blue!20,minimum size=4pt, inner sep=0pt, label=above:$y_3$] (y3) at (4, 5.5) {};
				\node[circle, draw, fill=blue!20,minimum size=4pt, inner sep=0pt, label=above:$y_2$] (y2) at (5.5, 5.5) {};
				\node[circle, draw, fill=blue!20,minimum size=4pt, inner sep=0pt, label=below:$y_1$] (y1) at (4, 4) {};
				\node[circle, draw, fill=blue!20,minimum size=4pt, inner sep=0pt, label=below:$y$] (y) at (5.5, 4) {};

				% Draw the edges
				\draw[thick, draw=blue] (x) -- (y);
				\draw[thick, draw=green] (x1) -- (x);
				\draw[thick, draw=red] (x) -- (x3);
				\draw[thick] (x) -- (x2);
				\draw[thick, draw=green] (y1) -- (y);
				\draw[thick, draw=green] (y) to (y2);
				\draw[thick] (y3) -- (y);
				\draw[thick, draw=green] (x1) to (y1);
				\draw[thick, bend right = 30, draw=red] (x3) to (y);
				\draw[thick, bend right = 30] (y3) to (x);

				\draw[thick] (-0.5, -0.) rectangle (6.5, 2.5);
				
				\node[circle, draw, fill=blue!20,minimum size=4pt, inner sep=0pt, label=above:$x^h$] (xh) at (0.2, 2) {};
				\node[circle, draw, fill=blue!20,minimum size=4pt, inner sep=0pt, label=below:$x_2^h$] (x2h) at (0.2, 1) {};
				\node[circle, draw, fill=blue!20,minimum size=4pt, inner sep=0pt, label=below:$y^h$] (yh) at (5.5, 1) {};
				\node[circle, draw, fill=blue!20,minimum size=4pt, inner sep=0pt, label=above:$y_2^h$] (y2h) at (5.5, 2) {};
				\node[circle, draw, fill=blue!20,minimum size=4pt, inner sep=0pt, label=below:$z$] (z) at (3, 1) {};
				\draw[decorate,draw=blue, decoration={snake, amplitude=0.02cm, segment length=0.1cm}] (xh) to (z);
				\draw[decorate,draw=red, decoration={snake, amplitude=0.02cm, segment length=0.1cm}] (yh) to (z);
				\draw[decorate, decoration={snake, amplitude=0.02cm, segment length=0.1cm}] (x2h) to (z);
				\draw[decorate,draw=green, decoration={snake, amplitude=0.02cm, segment length=0.1cm}] (y2h) to (z);
				\draw[thick, bend left=30, draw=blue] (xh) to (x);
				\draw[thick, bend right=30, draw=red] (yh) to (y);
				\draw[thick, bend left=30, draw=green] (y2) to (y2h);
				\draw[thick, bend right=30] (x2) to (x2h);
				\node at (3, -1.5) {\textbf{AQ$_3^1$}};
				\node at (3, 7) {\textbf{AQ$_3^0$}};
			\end{tikzpicture}
			\caption{Case 3.1}
			\label{fig:graph-example}
		\end{figure}
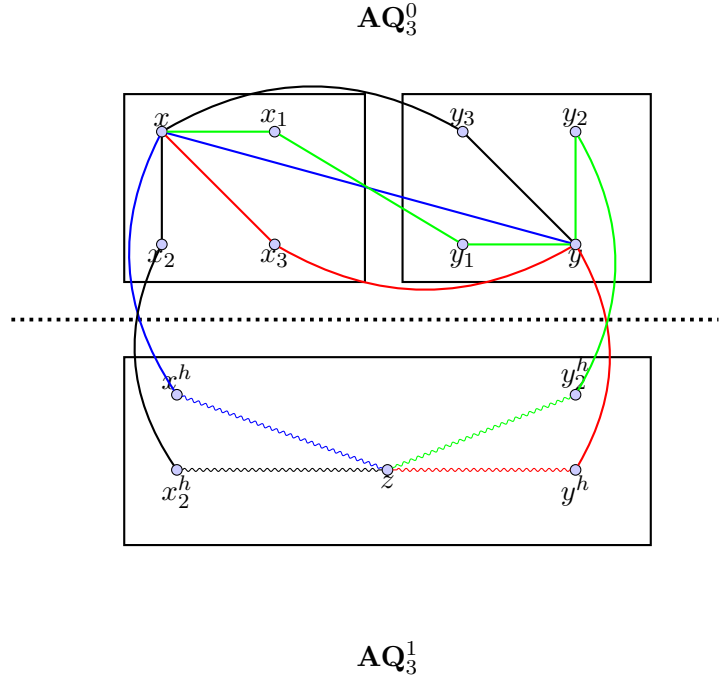
		\noindent
		$\psi_1=P_1\cup \{xx^h\}\cup P(x^h,z)\\
		\psi_2=P_2\cup \{yy^h\}\cup P(y^h,z)\\
		\psi_3=P_3\cup \{xx_2, x_2x_2^h\}\cup P(x_2^h,z)\\
		\psi_4=P_4\cup \{yy_2, y_2y_2^h\}\cup P(y_2^h,z) $\\
		
		\textbf{Subcase 3.2 :} Suppose  $x^{c_3}\neq y$. Then $x^h,y^h,x^c,y^c$ are all distinct vertices in $AQ_3^1.$ Also, there are $4-$vertex disjoint paths joining $x$ and $y$ in $AQ_3^0$, say, $P_1,P_2,P_3,$ and $P_4$. Let $X=\{x^h,y^h,x^c,y^c\}\subseteq V(AQ_3^1)$. Since $AQ_3^1$ is $4$-connected there exists four vertex disjoint paths $P(x^h,z),P(y^h,z),P(x^c,z)$, and $(y^c,z)$, each joining $z$ to every vertex of $X$.\\
		We get four vertex disjoint paths containing $x,y,$ and $z$ as follows, (See Figure 4).\\
		$\psi_1=P_1\cup \{xx^h\}\cup P(x^h,z)\\
		\psi_2=P_2\cup \{yy^h\}\cup P(y^h,z)\\
		\psi_3=P_3\cup \{xx^c\}\cup P(x_c,z)\\
		\psi_4=P_4\cup \{yy^c\} \cup P(y^c,z) $\\
		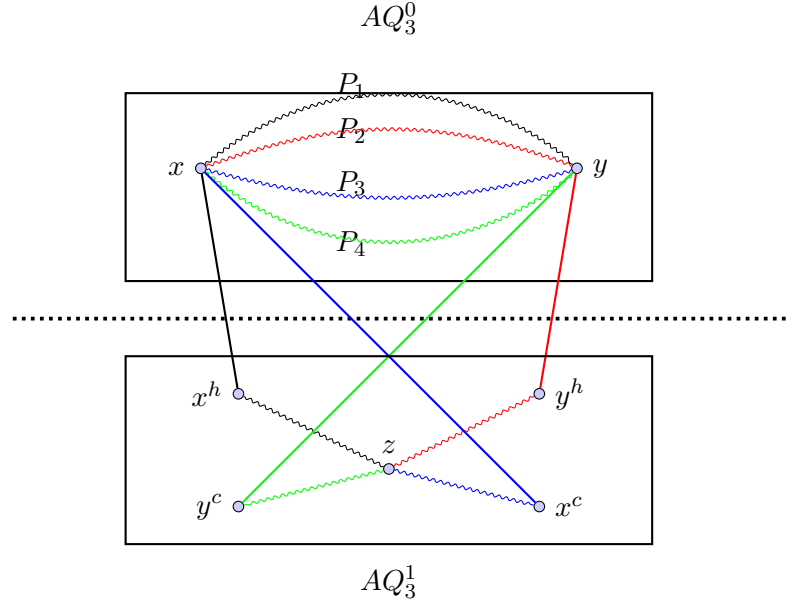
\begin{figure}[H]
			\centering
			\begin{tikzpicture}
				% Draw the first rectangle
				\draw[thick] (0.5, 0) rectangle (7.5, 2.5);
				% Label for the first rectangle
				\draw[dotted, line width=.5mm] (-1, -0.5) -- (9.4, -0.5);
				\node at (4, 3.5) {\textbf{$AQ_3^0$}};
				
				% Define nodes
				\node[circle, draw, fill=blue!20, minimum size=4pt, inner sep=0pt, label=right:$y$] (y) at (6.5, 1.5) {};
				\node[circle, draw, fill=blue!20, minimum size=4pt, inner sep=0pt, label=left:$x$] (x) at (1.5, 1.5) {};
				\node[circle, draw, fill=blue!20, minimum size=4pt, inner sep=0pt, label=above:$z$] (z) at (4, -2.5) {};
				\node[circle, draw, fill=blue!20, minimum size=4pt, inner sep=0pt, label=left:$x^h$] (xh) at (2, -1.5) {};
				\node[circle, draw, fill=blue!20, minimum size=4pt, inner sep=0pt, label=right:$x^c$] (xc) at (6, -3) {};
				\node[circle, draw, fill=blue!20, minimum size=4pt, inner sep=0pt, label=right:$y^h$] (yh) at (6, -1.5) {};
				\node[circle, draw, fill=blue!20, minimum size=4pt, inner sep=0pt, label=left:$y^c$] (yc) at (2, -3) {};
				
				% Labels for paths
				\node at (3.5, 2.6) {\textbf{$P_1$}};
				\node at (3.5, 2) {\textbf{$P_2$}};
				\node at (3.5, 1.3) {\textbf{$P_3$}};
				\node at (3.5, 0.5) {\textbf{$P_4$}};
				
				% Curly edges in the top rectangle
				\draw[decorate, decoration={snake, amplitude=0.02cm, segment length=0.1cm}, bend right=40, draw=green] (x) to (y);
				\draw[decorate, decoration={snake, amplitude=0.02cm, segment length=0.1cm}, bend right=15, draw=blue] (x) to (y);
				\draw[decorate, decoration={snake, amplitude=0.02cm, segment length=0.1cm}, bend left=20, draw=red] (x) to (y);
				\draw[decorate, decoration={snake, amplitude=0.02cm, segment length=0.1cm}, bend left=40] (x) to (y);
				\draw[decorate, decoration={snake, amplitude=0.02cm, segment length=0.1cm}, draw=red] (z) to (yh);
				\draw[decorate, decoration={snake, amplitude=0.02cm, segment length=0.1cm}] (z) to (xh);
				\draw[decorate, decoration={snake, amplitude=0.02cm, segment length=0.1cm}, draw=green] (z) to (yc);
				\draw[decorate, decoration={snake, amplitude=0.02cm, segment length=0.1cm}, draw=blue] (z) to (xc);
				\draw[thick] (xh) to (x);
				\draw[thick, draw=red] (y) to (yh);
				\draw[thick, draw=green] (y) to (yc);
				\draw[thick, draw=blue] (x) to (xc);
				
				% Draw the second rectangle below the first one
				\draw[thick] (0.5, -3.5) rectangle (7.5, -1);
				% Label for the second rectangle
				\node at (4, -4) {\textbf{$AQ_3^1$}};
			\end{tikzpicture}
			\caption{Case 3.2}
			\label{fig:graph-example}
		\end{figure}
		
	\end{proof}
	
	\begin{theorem}
		For $n\geq4$, $\pi_3(AQ_n)= \frac{3n}{2}-2$, when $n$ is even.
	\end{theorem}
	
	\begin{proof}
		By Lemma \ref{lemma6}, $\pi_3(AQ_n)\leq \frac{3n}{2}-2$, when $n$ is even. Now, we will prove that $\pi_3(AQ_n)= \frac{3n}{2}-2$ by constructing precisely  $\frac{3n}{2}-2 $ vertex disjoint paths containing any three distinct vertices in $V(AQ_n)$ for even $n\geq 4$.
		
		We prove this by using the mathematical induction on set of all even positive integers.\\
		For $n=4$, by Lemma \ref{lemma8}, $\pi _3(AQ_4)=\frac{3n}{2}-2=4$.\\
		Assume that this result is true for the augmented cube of dimension $n-2$, when $n$ is even.\\
		\textbf{Case 1 :} Consider, $\{x,y,z\}\subseteq V(AQ_{n-2}^{ij})$ for some $i$ and $j$, $0\leq i,j\leq 1$. By assumption there are $\frac{6(n-2)}{4}-2=\frac{3n}{2}-5$ vertex disjoint paths containing $x,y$ and $z$ in $AQ_{n-2}^{ij}$. Now, we have to find three additional paths containing $x,y$ and $z$ in $AQ_n$. Without loss of generality assume that $\{x,y,z\}\subseteq V(AQ_{n-2}^{00})$\\
		\textbf{Subcase 1.1 :} Suppose $x^{c_{n-2}}=y$ or $z$.\\
		Without loss of generality assume that $x^{c_{n-2}}=y$. Then $x^{c_{n-1}}=y^{h_{n-1}}$ and $x^{h_{n-1}}=y^{c_{n-1}}$.\\
		Observe that, for $\{x,y,z\}\subseteq V(AQ_{n-2}^{00})$, $\{x^h,y^h,z^h\}\subseteq V(AQ_{n-2}^{10})$, $\{x^c,y^c,z^c\}\subseteq V(AQ_{n-2}^{11})$ and $\{x^{h_{n-1}},y^{h_{n-1}},z^{h_{n-1}},z^{c_{n-1}}\}\subseteq V(AQ_{n-2}^{01})$. $z^c$ has adjecent vertex $a \in V(AQ_{n-2}^{10})$ which is not adjecent to $x^h,y^h$ or $z^h$. Then, by Lemma \ref{lemma7}, we get two vertex disjiont paths $P(x^h,z^h)$ and $P(y^h,a)$ in $AQ_{n-2}^{10}$. As  $AQ_{n-2}^{01}\diamond AQ_{n-2}^{11}$ is $2n-4$ connected, by Lemma \ref{lemma7}, there are two vertex disjoint paths $P(x^c,z^{c_{n-1}})$ and $P(y^c,z^{h_{n-1}})$ in $AQ_{n-2}^{01}\diamond AQ_{n-2}^{11}$ which does not conatain $x^{h_{n-1}},y^{h_{n-1}}$, and $z^c$ as there internal vertices.\\
		By induction, we have already $\frac{3n}{2}-5$ vertex disjoint paths containing $x,y,$ and $z$. Remainning three vertex disjoint paths containing $x,y$ and $z$ in $AQ_n$ are as follows, (See Figure 5).\\
		$
		\psi_{\frac{3n}{2}-4} = \{yx^{h_{n-1}}, xx^{h_{n-1}}, xx^c\} \cup P(x^c, z^{c_{n-1}})\cup \{zz^{c_{n-1}}\}, \\ 
		\psi_{\frac{3n}{2}-3} = \{xy^{h_{n-1}}, yy^{h_{n-1}}, yy^c\} \cup P(y^c, z^{h_{n-1}})\cup \{zz^{h_{n-1}}\}  , \\ 
		\psi_{\frac{3n}{2}-2} = \{xx^h\}\cup P(x^h, z^h) \cup \{z^hz, z^cz, az^c\} \cup  P(y^h, a) \cup \{yy^h\}.
		$
		\\

		\begin{figure}[H]
			\centering
			\begin{tikzpicture}
				% Draw the top-left rectangle
				\draw[thick] (-7, 2) rectangle (-0.5, 5);
				
				% Draw the top-right rectangle
				\draw[thick] (-0.2, 2) rectangle (6.2, 5);
				
				% Draw the bottom-left rectangle
				\draw[thick] (-7, -2.5) rectangle (-0.5, 0.5);
				
				% Draw the bottom-right rectangle
				\draw[thick] (-0.2, -2.5) rectangle (6.2, 0.5);
				\draw[dotted, line width=.5mm] (-7.5, 1.3) -- (7., 1.3);
				% Top-left rectangle nodes
				\node[circle, draw, fill=black!20, minimum size=4pt, inner sep=0pt, label={[font=\scriptsize]left:$x$}] (x) at (-6, 4.4) {};
				\node[circle, draw, fill=black!20, minimum size=4pt, inner sep=0pt, label={[font=\scriptsize]left:$y$}] (y) at (-1.5, 2.5) {};
				\node[circle, draw, fill=black!20, minimum size=4pt, inner sep=0pt, label={[font=\scriptsize]left:$z$}] (z) at (-3.5, 3.5) {};
				
				% Top-right rectangle nodes
				\node[circle, draw, fill=black!20, minimum size=4pt, inner sep=0pt, label={[font=\scriptsize]left:$x^{h_{n-1}}$}] (xhn-1) at (1, 4.4) {};
				\node[circle, draw, fill=black!20, minimum size=4pt, inner sep=0pt, label={[font=\scriptsize]right:$y^{h_{n-1}}$}] (yhn-1) at (5.5, 2.5) {};
				\node[circle, draw, fill=black!20, minimum size=4pt, inner sep=0pt, label={[font=\scriptsize]left:$z^{h_{n-1}}$}] (zhn-1) at (2, 2.5) {};
				\node[circle, draw, fill=black!20, minimum size=4pt, inner sep=0pt, label={[font=\scriptsize]right:$z^{c_{n-1}}$}] (zcn-1) at (4.5, 4) {};
				
				% Bottom-left rectangle nodes
				\node[circle, draw, fill=black!20, minimum size=4pt, inner sep=0pt, label={[font=\scriptsize]left:$x^h$}] (xh) at (-6, -0.1) {};
				\node[circle, draw, fill=black!20, minimum size=4pt, inner sep=0pt, label={[font=\scriptsize]right:$y^h$}] (yh) at (-1.5, -1.9) {};
				\node[circle, draw, fill=black!20, minimum size=4pt, inner sep=0pt, label={[font=\scriptsize]left:$z^h$}] (zh) at (-3.5, -1) {};
				\node[circle, draw, fill=black!20, minimum size=4pt, inner sep=0pt, label={[font=\scriptsize]left:$a$}] (a) at (-5, -1.5) {};
				
				% Bottom-right rectangle nodes
				\node[circle, draw, fill=black!20, minimum size=4pt, inner sep=0pt, label={[font=\scriptsize]left:$y^c$}] (yc) at (1, -0.1) {};
				\node[circle, draw, fill=black!20, minimum size=4pt, inner sep=0pt, label={[font=\scriptsize]left:$x^c$}] (xc) at (5.5, -1.9) {};
				\node[circle, draw, fill=black!20, minimum size=4pt, inner sep=0pt, label={[font=\scriptsize]right:$z^c$}] (zc) at (3.5, -1) {};
				
				% Draw edges
				\draw[thin,red,bend left=15] (x) to (xhn-1);
				\draw[thin,red] (y) to (xhn-1);
				\draw[thin,green] (x) to (yhn-1);
				\draw[thin,green,bend right=15] (y) to (yhn-1);
				
				\draw[decorate, decoration={snake, amplitude=0.02cm, segment length=0.1cm},red,bend left=30] (xc) to (zcn-1);
				\draw[decorate, decoration={snake, amplitude=0.02cm, segment length=0.1cm},green] (yc) to (zhn-1);
				\draw[thin,red] (z) to (zcn-1);
				\draw[thin,green] (z) to (zhn-1);
				\draw[decorate, decoration={snake, amplitude=0.02cm, segment length=0.1cm},blue] (xh) to (zh);
				\draw[decorate, decoration={snake, amplitude=0.02cm, segment length=0.1cm},blue] (yh) to (a);
				\draw[thin,blue] (zc) to (a);
				% Labels for each rectangle
				\node at (-4, 5.5) {\textbf{$AQ_{n-2}^{00}$}};
				\node at (3, 5.5) {\textbf{$AQ_{n-2}^{01}$}};
				\node at (-4, -3) {\textbf{$AQ_{n-2}^{10}$}};
				\node at (3, -3) {\textbf{$AQ_{n-2}^{11}$}};
				\draw[thin,blue] (x) to (xh);
				\draw[thin,blue] (y) to (yh);
				\draw[thin,green] (y) to (yc);
				\draw[thin, red, bend right=50] (x) to (xc);
				\draw[thin,blue] (z) to (zh);
				\draw[thin,blue] (z) to (zc);
			\end{tikzpicture}
			\caption{Case 1.1}
			\label{fig:case-1.1}
		\end{figure}
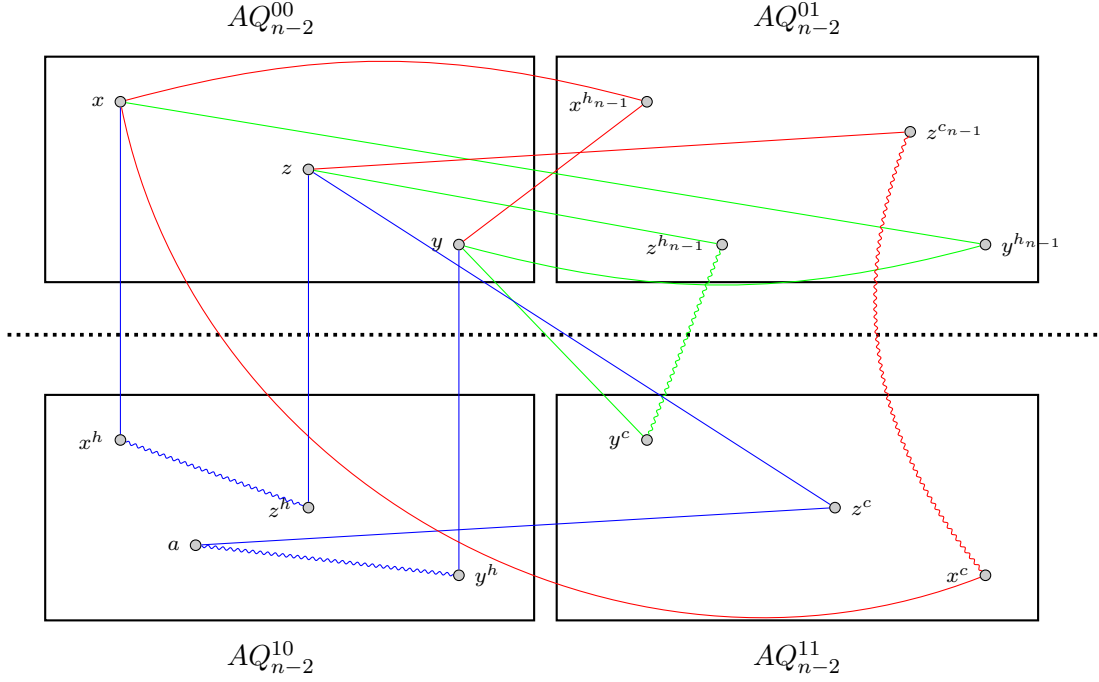
		
		\vspace*{-10pt}

		\textbf{Subcase 1.2 :} Suppose $x^{c_{n-2}}\neq y$ and $z$.\\
		Then  $x^{h_{n-1}},y^{h_{n-1}},z^{h_{n-1}},x^{c_{n-1}},y^{c_{n-1}},z^{c_{n-1}}$ are distinct vertices of $AQ_{n-2}^{01}.$  There exists three vertex disjoint paths $P(x^{h_{n-1}},z^{h_{n-1}}),P(y^{h_{n-1}},z^{c_{n-1}}) $ and $P(x^{c_{n-1}},y^{c_{n-1}})$ in $AQ_{n-2}^{01}.$ Since, $\{x,y,z\}\subseteq V(AQ_{n-2}^{00}),$ $\{x^h,y^h,z^h\}\subseteq V(AQ_{n-2}^{10})$ and $\{x^c,y^c,z^c\}\subseteq V(AQ_{n-2}^{11})$. By Lemma \ref{lemma7}, we get three vertex disjoint paths $P(x^h,z^h),P(y^h,z^c)$ and $P(x^c,y^c)$ in $AQ_{n-1}^1.$\\
		By induction, we have already $\frac{3n}{2}-2$ vertex disjoint paths containing $x,y,$ and $z$. Remainning three vertex disjoint paths containing $x,y$ and $z$ in $AQ_n$ as follows, (See figure 6).\\
		$
		\psi_{\frac{3n}{2}-4} = \{xx^{h_{n-1}}\}\cup P(x^{h_{n-1}}, z^{h_{n-1}})  \cup \{zz^{h_{n-1}}, zz^{c_{n-1}}\}  \cup P(y^{h_{n-1}}, z^{c_{n-1}})\cup \{yy^{h_{n-1}}\}, \\ 
		\psi_{\frac{3n}{2}-3} =\{yy^{c_{n-1}}\}\cup P(x^{c_{n-1}}, y^{c_{n-1}}) \cup  \{xx^{c_{n-1}}, xx^h\}\cup P(x^h, z^h) \cup \{zz^h\} , \\ 
		\psi_{\frac{3n}{2}-2} = \{xx^c\}\cup P(x^c, y^c)  \cup \{yy^c, y^hy\} \cup P(y^h, z^c)  \cup \{z^cz\}.
		$
		\\
		\begin{figure}[H]
			\centering
			\begin{tikzpicture}
				% Draw the top-left rectangle
				\draw[thick] (-7, 2) rectangle (-0.5, 5);
				
				% Draw the top-right rectangle
				\draw[thick] (-0.2, 2) rectangle (6.2, 5);
				
				% Draw the bottom-left rectangle
				\draw[thick] (-7, -2.5) rectangle (-0.5, 0.5);
				
				% Draw the bottom-right rectangle
				\draw[thick] (-0.2, -2.5) rectangle (6.2, 0.5);
				\draw[dotted, line width=.5mm] (-7.4, 1.3) -- (7, 1.3);
				% Top-left rectangle nodes
				\node[circle, draw, fill=black!20, minimum size=4pt, inner sep=0pt, label={[font=\scriptsize]left:$x$}] (x) at (-5, 4.4) {};
				\node[circle, draw, fill=black!20, minimum size=4pt, inner sep=0pt, label={[font=\scriptsize]left:$y$}] (y) at (-2.5, 2.5) {};
				\node[circle, draw, fill=black!20, minimum size=4pt, inner sep=0pt, label={[font=\scriptsize]left:$z$}] (z) at (-3.5, 3.5) {};
				
				% Top-right rectangle nodes
				\node[circle, draw, fill=black!20, minimum size=4pt, inner sep=0pt, label={[font=\scriptsize]left:$x^{h_{n-1}}$}] (xhn-1) at (1, 4.4) {};
				\node[circle, draw, fill=black!20, minimum size=4pt, inner sep=0pt, label={[font=\scriptsize]right:$x^{c_{n-1}}$}] (xcn-1) at (5.5, 2.5) {};
				\node[circle, draw, fill=black!20, minimum size=4pt, inner sep=0pt, label={[font=\scriptsize]left:$z^{h_{n-1}}$}] (zhn-1) at (2, 2.5) {};
				\node[circle, draw, fill=black!20, minimum size=4pt, inner sep=0pt, label={[font=\scriptsize]right:$z^{c_{n-1}}$}] (zcn-1) at (4.5, 4.4) {};
				\node[circle, draw, fill=black!20, minimum size=4pt, inner sep=0pt, label={[font=\scriptsize]right:$y^{c_{n-1}}$}] (ycn-1) at (2.5, 4.4) {};
				\node[circle, draw, fill=black!20, minimum size=4pt, inner sep=0pt, label={[font=\scriptsize]right:$y^{h_{n-1}}$}] (yhn-1) at (3.5, 2.5) {};
				
				% Bottom-left rectangle nodes
				\node[circle, draw, fill=black!20, minimum size=4pt, inner sep=0pt, label={[font=\scriptsize]left:$x^h$}] (xh) at (-6, -0.1) {};
				\node[circle, draw, fill=black!20, minimum size=4pt, inner sep=0pt, label={[font=\scriptsize]right:$y^h$}] (yh) at (-1.5, -1.9) {};
				\node[circle, draw, fill=black!20, minimum size=4pt, inner sep=0pt, label={[font=\scriptsize]left:$z^h$}] (zh) at (-3.5, -1) {};
				
				% Bottom-right rectangle nodes
				\node[circle, draw, fill=black!20, minimum size=4pt, inner sep=0pt, label={[font=\scriptsize]left:$y^c$}] (yc) at (1, -0.1) {};
				\node[circle, draw, fill=black!20, minimum size=4pt, inner sep=0pt, label={[font=\scriptsize]left:$x^c$}] (xc) at (5.5, -1.9) {};
				\node[circle, draw, fill=black!20, minimum size=4pt, inner sep=0pt, label={[font=\scriptsize]right:$z^c$}] (zc) at (3.5, -1) {};
				
				% Draw edges
				\draw[decorate, decoration={snake, amplitude=0.02cm, segment length=0.1cm},green] (xhn-1) to (zhn-1);
				\draw[decorate, decoration={snake, amplitude=0.02cm, segment length=0.1cm},green] (yhn-1) to (zcn-1);
				\draw[thin,green,bend right=15] (xhn-1) to (x);
				\draw[thin,green,bend left=15] (yhn-1) to (y);
				\draw[thin,green] (zhn-1) to (z);
				\draw[thin,green] (zcn-1) to (z);
				\draw[decorate, decoration={snake, amplitude=0.02cm, segment length=0.1cm},red] (xh) to (zh);
				\draw[decorate, decoration={snake, amplitude=0.02cm, segment length=0.1cm},red] (xcn-1) to (ycn-1);
				\draw[thin,red] (xcn-1) to (x);
				\draw[thin,red] (ycn-1) to (y);
				\draw[decorate, decoration={snake, amplitude=0.02cm, segment length=0.1cm},blue,bend right=25] (xc) to (yc);
				\draw[decorate, decoration={snake, amplitude=0.02cm, segment length=0.1cm},blue] (zc) to (yh);
				% Labels for each rectangle
				\node at (-4, 5.5) {\textbf{$AQ_{n-2}^{00}$}};
				\node at (3, 5.5) {\textbf{$AQ_{n-2}^{01}$}};
				\node at (-4, -3) {\textbf{$AQ_{n-2}^{10}$}};
				\node at (3, -3) {\textbf{$AQ_{n-2}^{11}$}};
				\draw[thin,red] (x) to (xh);
				\draw[thin,blue] (y) to (yh);
				\draw[thin,blue] (y) to (yc);
				\draw[thin, blue, bend right=50] (x) to (xc);
				\draw[thin,red] (z) to (zh);
				\draw[thin,blue] (z) to (zc);
			\end{tikzpicture}
			\caption{Case 1.2}
			\label{fig:case-1.2}
		\end{figure}
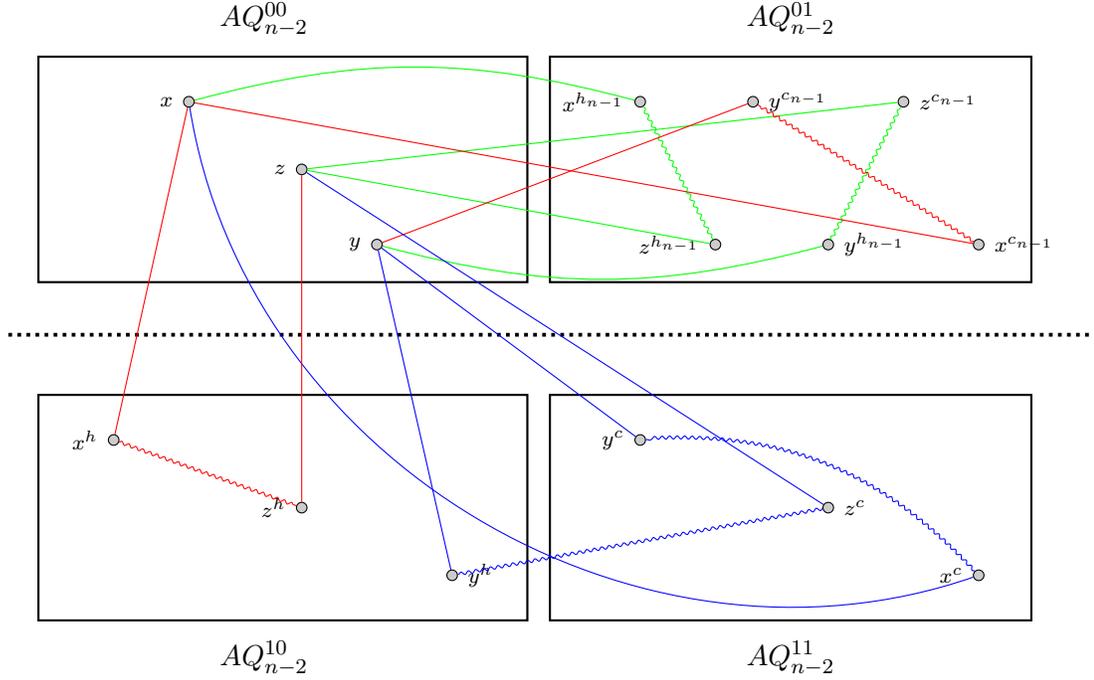

		\textbf{Case 2 : } $x,y\in V(AQ_{n-2}^{00})$ and $z\in V(AQ_{n-2}^{01}) $\\
		\textbf{Subcase 2.1 :} Suppose, $z^{c_{n-1}}=x$ or $y$.\\
		Without loss of generality consider, $z^{c_{n-1}}=x$ which implies that $z^c=x^h$ and $z^h=x^c.$\\
		As, $AQ_{n-2}^{00}$ is $2n-5$ connected, for, $x,y\in V(AQ_{n-2}^{00})$ there exists $2n-5$ vertex disjoint paths joining $x$ and $y$ in $AQ_{n-2}^{00},$ say, $Q_1,Q_2,\dots ,Q_{n-3},P_1,P_2,\dots ,P_{n-2}.$
		Let $x_1,x_2,\dots ,x_{n-3}$ be the neighbours of x along $Q_1,Q_2,\dots ,Q_{n-3}$ respectively and $y_1,y_2,\dots ,y_{n-3}$ be the neighbours of y along $Q_1,Q_2,\dots ,Q_{n-3}$ respectively. Note that $\{x_1^{h_{n-1}},x_2^{h_{n-1}},\dots ,x_{n-3}^{h_{n-1}}\}$ is the set of neighbours of $x^{h_{n-1}}$ in $AQ_{n-2}^{01}$ and $\{y_1^{h_{n-1}},y_2^{h_{n-1}},\dots ,y_{n-3}^{h_{n-1}}\}$ is the set of neighbours of $y^{h_{n-1}}$ in $AQ_{n-2}^{01}$. As $AQ_{n-2}^{00}$ is $2n-5$ connected, by Lemma \ref{lemma1}, there exist $2n-5$ vertex disjoint paths joining $z$ to every vertex in $\{x_1^{h_{n-1}},x_2^{h_{n-1}},\dots ,x_{n-3}^{h_{n-1}},y_1^{h_{n-1}},y_2^{h_{n-1}},\dots ,y_{n-3}^{h_{n-1}},y^{h_{n-1}} \}$, say, $P(x_i^{h_{n-1}},z),i=1,2,\dots ,n-3$, $P(y_i^{h_{n-1}},z), i=1,2,\dots ,n-3 $ and $P(y^{h_{n-1}},z)$.\\
		Then we get $2n-\frac{n}{2}-2$ vertex disjoint paths containing $x,y$ and $z$ in $AQ_n$ as follows, (See Figure 7).\\
		$
		\psi_i = P_i\cup \{xx_i\} \cup \{x_ix_i^{h_{n-1}}\} \cup P(x_i^{h_{n-1}}, z), 1\leq i  \leq \frac{n}{2}-2  \\
		\psi_{\frac{n}{2}-1} = P_{\frac{n}{2}-1}\cup \{yy_1\} \cup \{y_1y_1^{h_{n-1}}\} \cup P(y_1^{h_{n-1}}, z), \\ 
		\vdots \\ 
		\psi_{n-4} = P_{n-4}\cup \{yy_{\frac{n}{2}-2} , y_{\frac{n}{2}-2}y_{\frac{n}{2}-2}^{h_{n-1}}\} \cup P(y_{\frac{n}{2}-2}^{h_{n-1}}, z), \\ 
		\psi_{n-3} = \{x_{\frac{n}{2}-1}x_{\frac{n}{2}-1}^{h_{n-1}}\} \cup P(z, x_{\frac{n}{2}-1}^{h_{n-1}}) \cup \{y_{\frac{n}{2}-1}y_{\frac{n}{2}-1}^{h_{n-1}}\} \cup P(y_{\frac{n}{2}-1}^{h_{n-1}}, z), \\ 
		\vdots \\ 
		\psi_{n+\frac{n}{2}-5} = \{xx_{n-3}\} \cup \{x_{n-3}x_{n-3}^{h_{n-1}}\} \cup P(z, x_{n-3}^{h_{n-1}}) \cup \{y_{n-3}y_{n-3}^{h_{n-1}}\} \cup P(y_{n-3}^{h_{n-1}}, z) \cup  \{yy_{n-3}\}, \\ 
		\psi_{n+\frac{n}{2}-4} = \{yy^{h_{n-1}}\} \cup P(z, y^{h_{n-1}}) \cup \{zx\}, \\ 
		\psi_{n+\frac{n}{2}-3} = P_{n-3}\cup \{xx^h\} \cup \{zz^c\}  , \\ 
		\psi_{n+\frac{n}{2}-2} = P_{n-2}\cup \{xx^c\} \cup \{zz^h\}  .
		$
		\\
		\vspace{-10pt}
		\begin{figure}[H]
			\centering
			
			\begin{tikzpicture}{baseline}
				% Draw the top-left rectangle
				\draw[thick] (-7, 2) rectangle (-0.5, 5);
				
				% Draw the top-right rectangle
				\draw[thick] (-0.2, 2) rectangle (6.2, 5);
				
				% Draw the bottom-left rectangle
				\draw[thick] (-7, -2.5) rectangle (-0.5, 0.5);
				
				% Draw the bottom-right rectangle
				\draw[thick] (-0.2, -2.5) rectangle (6.2, 0.5);
				
				% Labels for each rectangle
				\node at (-4, 5.8) {\textbf{$AQ_{n-2}^{00}$}};
				\node at (3, 5.8) {\textbf{$AQ_{n-2}^{01}$}};
				\node at (-4, -3) {\textbf{$AQ_{n-2}^{10}$}};
				\node at (3, -3) {\textbf{$AQ_{n-2}^{11}$}};
				
				% Dotted line
				\draw[dotted, line width=.5mm] (-8, 1.3) -- (7, 1.3);
				
				\node[circle, draw, fill=black!20, minimum size=4pt, inner sep=0pt, label=left:$x$] (x) at (-6.4, 4.6) {};
				\node[circle, draw, fill=black!20, minimum size=4pt, inner sep=0pt, label=above:$x^h(z^c)$] (xhzc) at (-6.4, -0.2) {};
				\node[circle, draw, fill=black!20, minimum size=4pt, inner sep=0pt, label=above:$y^c$] (yc) at (0.3, -0.3) {};
				\node[circle, draw, fill=black!20, minimum size=2pt, inner sep=0pt, label={[font=\scriptsize]right:$x_1$}] (x1) at (-5.2, 4.8) {};
				\node[circle, draw, fill=black!20, minimum size=2pt, inner sep=0pt, label={[font=\scriptsize]right:$x_2$}] (x2) at (-5.2, 4.5) {};
				\node[circle, draw, fill=black!20, minimum size=1pt, inner sep=0pt] (a) at (-5.2, 4.3) {};
				\node[circle, draw, fill=black!20, minimum size=1pt, inner sep=0pt] (a) at (-5.2, 4.4) {};
				\node[circle, draw, fill=black!20, minimum size=1pt, inner sep=0pt] (a) at (-5.2, 4.1) {};
				\node[circle, draw, fill=black!20, minimum size=1pt, inner sep=0pt] (a) at (-5.2, 4) {};
				\node[circle, draw, fill=black!20, minimum size=2pt, inner sep=0pt, label={[font=\scriptsize]right:$x_{n-3}$}] (xn-3) at (-5.2, 3.9) {};
				\node[circle, draw, fill=black!20, minimum size=2pt, inner sep=0pt, label={[font=\scriptsize]right:$x_{\frac{n}{2}-2}$}] (xn/2-2) at (-5.2, 4.2) {};
				\node[circle, draw, fill=black!20, minimum size=2pt, inner sep=0pt, label={[font=\scriptsize]above:$y_1$}] (y1) at (-2.2, 4.4) {};
				\node[circle, draw, fill=black!20, minimum size=2pt, inner sep=0pt, label={[font=\scriptsize]left:$y_2$}] (y2) at (-2.3, 4) {};
				\node[circle, draw, fill=black!20, minimum size=1pt, inner sep=0pt] (a) at (-2.34, 3.8) {};
				\node[circle, draw, fill=black!20, minimum size=1pt, inner sep=0pt] (a) at (-2.36, 3.7) {};
				\node[circle, draw, fill=black!20, minimum size=2pt, inner sep=0pt, label={[font=\scriptsize]left:$y_{n-3}$}] (yn-3) at (-2.5, 3.3) {};
				\node[circle, draw, fill=black!20, minimum size=1pt, inner sep=0pt] (a) at (-2.46, 3.4) {};
				\node[circle, draw, fill=black!20, minimum size=1pt, inner sep=0pt] (a) at (-2.44, 3.5) {};
				\node[circle, draw, fill=black!20, minimum size=2pt, inner sep=0pt, label={[font=\scriptsize]left:$y_{\frac{n}{2}-2}$}] (yn/2-2) at (-2.4, 3.6) {};
				\node[circle, draw, fill=blue!20, minimum size=4pt, inner sep=0pt, label=below:$y$] (y) at (-1.1, 3.6) {};
				\node[circle, draw, fill=black!20, minimum size=4pt, inner sep=0pt, label=below:$y^h$] (yh) at (-1.1, -1.6) {};
				\node[circle, draw, fill=black!20, minimum size=4pt, inner sep=0pt, label=below:$x^c(z^h)$] (xc) at (5, -1.6) {};
				\draw[decorate, decoration={snake, amplitude=0.02cm, segment length=0.1cm},blue] (xhzc) to (yh);
				\draw[decorate, decoration={snake, amplitude=0.02cm, segment length=0.1cm},yellow] (yc) to (xc);
				\draw[thin,yellow]  (y)to(yc);
				\draw[thin,red] (x) to (x1);
				\draw[thin,red] (x) to (x2);
				\draw[thin,green] (x) to (xn-3);
				\draw[thin] (y) to (y1);
				\draw[thin] (y) to (y2);
				\draw[thin,green] (y) to (yn-3);
				\draw[thin] (y) to (yn/2-2);
				\draw[dotted] (x1) to (y1);
				\draw[dotted] (x2) to (y2);
				\draw[dotted] (yn-3) to (xn-3);
				\draw[decorate, decoration={snake, amplitude=0.02cm, segment length=0.1cm}, bend left=35,red] (y) to (x);
				\draw[decorate, decoration={snake, amplitude=0.02cm, segment length=0.1cm}, bend left=50] (y) to (x);
				\draw[thin] (x) to (xn/2-2);
				\draw[dotted] (yn/2-2) to (xn/2-2);
				\ draw[thin, bend right=8] (yn/2-2) to (y);
				\draw[decorate, decoration={snake, amplitude=0.02cm, segment length=0.1cm}, bend left=90,yellow] (y) to (x);
				\draw[decorate, decoration={snake, amplitude=0.02cm, segment length=0.1cm}, bend left=70,green] (y) to (x);
				\node at (-4.7, 3.5) {{\scriptsize {$P_1$}}};
				\node at (-4.7, 3.1) {\scriptsize {$P_2$}};
				\node at (-4.6, 2.5) {\scriptsize{$P_{n-2}$}};
				\node at (-3.8, 4.7) {\scriptsize{$Q_1$}} ;
				\node at (-3.8, 4.3) {\scriptsize{$Q_2$}};
				\node at (-3.8, 3.6) {\scriptsize{$Q_{n-3}$}};
				\node[circle, draw, fill=black!20, minimum size=1pt, inner sep=0pt] (a) at (-3.8, 4.0) {};
				\node[circle, draw, fill=black!20, minimum size=1pt, inner sep=0pt] (a) at (-3.8, 4.1) {};
				\node[circle, draw, fill=black!20, minimum size=1pt, inner sep=0pt] (a) at (-3.8, 3.9) {};
				\node[circle, draw, fill=black!20, minimum size=1pt, inner sep=0pt] (a) at (-4.7, 2.9) {};
				\node[circle, draw, fill=black!20, minimum size=1pt, inner sep=0pt] (a) at (-4.7, 2.8) {};
				\node[circle, draw, fill=black!20, minimum size=1pt, inner sep=0pt] (a) at (-4.7, 2.7) {};
				\node[circle, draw, fill=black!20, minimum size=1pt, inner sep=0pt] (a) at (-4.7, 2.6) {};
				
				\node[circle, draw, fill=black!20, minimum size=2pt, inner sep=0pt, 
				label={[font=\scriptsize, align=left]left:$x_1^{h_{n-1}}$}] (x1hn-1) at (1, 4.5) {};
				\node[circle, draw, fill=black!20, minimum size=1pt, inner sep=0pt] (a) at (1, 4.2) {};
				\node[circle, draw, fill=black!20, minimum size=1pt, inner sep=0pt] (a) at (1, 4.1) {};
				\node[circle, draw, fill=black!20, minimum size=1pt, inner sep=0pt] (a) at (1, 4) {};
				\node[circle, draw, fill=black!20, minimum size=1pt, inner sep=0pt] (a) at (1, 4.3) {};
				\node[circle, draw, fill=black!20, minimum size=2pt, inner sep=0pt, 
				label={[font=\scriptsize, align=left]left:$x_{\frac{n}{2}-2}^{h_{n-1}}$}] (xn/2-2hn-1) at (1,3.9) {};
				\node[circle, draw, fill=black!20, minimum size=2pt, inner sep=0pt, 
				label={[font=\scriptsize, align=left]left:$x_{\frac{n}{2}-1}^{h_{n-1}}$}] (xn/2-1hn-1) at (1,3.3) {};
				\node[circle, draw, fill=black!20, minimum size=1pt, inner sep=0pt] (a) at (1, 3.2) {};
				\node[circle, draw, fill=black!20, minimum size=1pt, inner sep=0pt] (a) at (1, 3.1) {};
				\node[circle, draw, fill=black!20, minimum size=1pt, inner sep=0pt] (a) at (1, 3) {};
				\node[circle, draw, fill=black!20, minimum size=2pt, inner sep=0pt, 
				label={[font=\scriptsize, align=left]left:$x_{n-3}^{h_{n-1}}$}] (xn-3hn-1) at (1,2.7) {};

				\node[circle, draw, fill=black!20, minimum size=2pt, inner sep=0pt, 
				label={[font=\scriptsize, align=left]left:$y_1^{h_{n-1}}$}] (y1hn-1) at (4.5, 4.5) {};
				\node[circle, draw, fill=black!20, minimum size=1pt, inner sep=0pt] (a) at (4.5, 4.2) {};
				\node[circle, draw, fill=black!20, minimum size=1pt, inner sep=0pt] (a) at (4.5, 4.1) {};
				\node[circle, draw, fill=black!20, minimum size=1pt, inner sep=0pt] (a) at (4.5, 4) {};
				\node[circle, draw, fill=black!20, minimum size=1pt, inner sep=0pt] (a) at (4.5, 4.3) {};
				\node[circle, draw, fill=black!20, minimum size=2pt, inner sep=0pt, 
				label={[font=\scriptsize, align=left]left:$y_{\frac{n}{2}-2}^{h_{n-1}}$}] (yn/2-2hn-1) at (4.5,3.9) {};
				\node[circle, draw, fill=black!20, minimum size=2pt, inner sep=0pt, 
				label={[font=\scriptsize, align=right]left:$y_{\frac{n}{2}-1}^{h_{n-1}}$}] (yn/2-1hn-1) at (4.5,3.3) {};
				\node[circle, draw, fill=black!20, minimum size=1pt, inner sep=0pt] (a) at (4.5, 3.2) {};
				\node[circle, draw, fill=black!20, minimum size=1pt, inner sep=0pt] (a) at (4.5, 3.1) {};
				\node[circle, draw, fill=black!20, minimum size=1pt, inner sep=0pt] (a) at (4.5, 3) {};
				\node[circle, draw, fill=black!20, minimum size=2pt, inner sep=0pt, 
				label={[font=\scriptsize, align=left]left:$y_{n-3}^{h_{n-1}}$}] (yn-3hn-1) at (4.5,2.7) {};
				\node[circle, draw, fill=black!20, minimum size=4pt, inner sep=0pt, 
				label=right:$z$] (z) at (5.5,2.37) {};
				\node[circle, draw, fill=black!20, minimum size=2pt, inner sep=0pt, 
				label={[font=\scriptsize, align=right]right:$y^{h_{n-1}}$}] (yhn-1) at (5,3.8) {};
				
				\draw[decorate, decoration={snake, amplitude=0.02cm, segment length=0.1cm}, thin, green] (z) -- (yn-3hn-1);
				\draw[decorate, decoration={snake, amplitude=0.02cm, segment length=0.1cm}, thin,green] (z) -- (yn/2-1hn-1);
				\draw[decorate, decoration={snake, amplitude=0.02cm, segment length=0.1cm}, thin] (z) -- (y1hn-1);
				\draw[decorate, decoration={snake, amplitude=0.02cm, segment length=0.1cm}, thin,] (z) -- (yn/2-2hn-1);
				\draw[decorate, decoration={snake, amplitude=0.02cm, segment length=0.1cm}, thin, bend left=20,green] (z) to (xn-3hn-1);
				\draw[decorate, decoration={snake, amplitude=0.02cm, segment length=0.1cm},green, thin, bend left=20] (z) to (xn/2-1hn-1);
				\draw[decorate, decoration={snake, amplitude=0.02cm, segment length=0.1cm}, thin, bend left=20,red] (z) to (xn/2-2hn-1);
				\draw[decorate, decoration={snake, amplitude=0.02cm, segment length=0.1cm}, thin, bend left=20, red] (z) to (x1hn-1);
				\draw[decorate, decoration={snake, amplitude=0.02cm, segment length=0.1cm}, thin,orange] (z) -- (yhn-1);
				
				\draw[thin,bend left=45,orange] (z)to(x);
				\draw[thin,bend right=25,red] (x1hn-1)to(x1);
				\draw[thin,bend right=25,green]  (yn-3)to(yn-3hn-1);
				
				\draw[thin,bend right=25,green]  (xn-3)to(xn-3hn-1);
				\draw[thin,yellow]  (z)to(xc);
				\draw[thin,blue]  (x)to(xhzc);
				\draw[thin,blue]  (y)to(yh);
				\draw[thin,bend left=30]  (y1)to(y1hn-1);
			\end{tikzpicture}
			
			\caption{Case 2.1}
			\label{fig:enter-label}
		\end{figure}
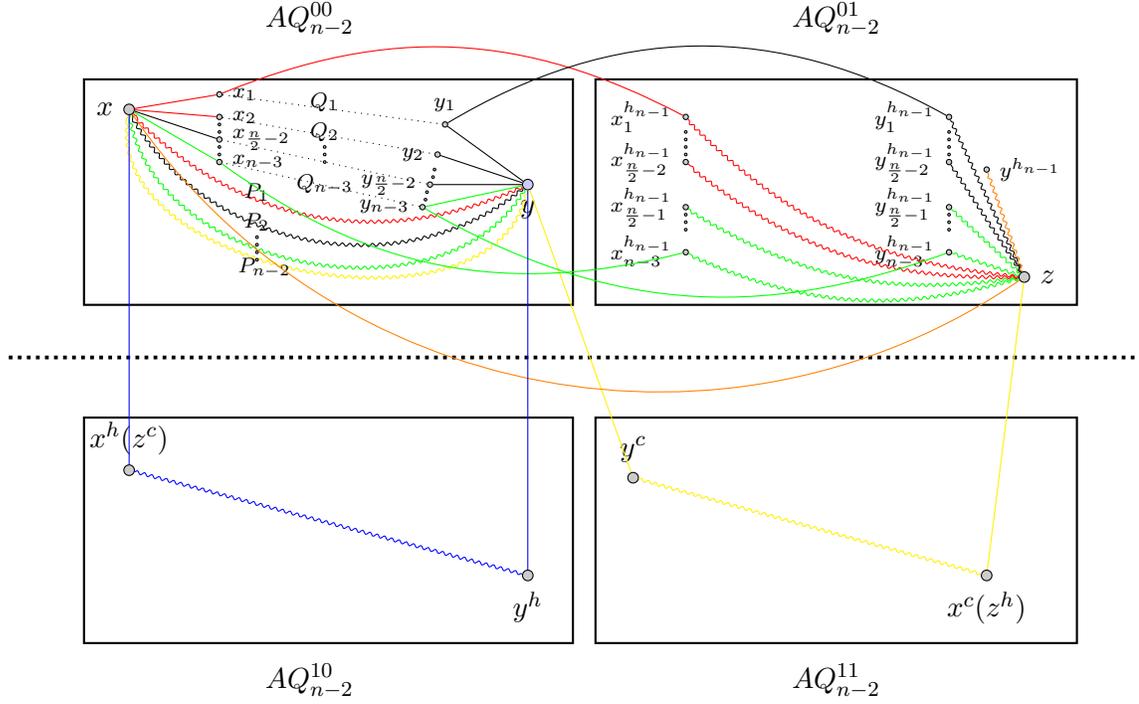  
		\vspace{-10pt}
		
		\noindent
		\textbf{Subcase 2.2 :} $z^{c_{n-1}}\neq x $ and $y$, which implies that $z^h, x^h,y^h$ and $z^c, x^c,y^c$ are distinct vertices in $AQ_{n-1}^1$. By Lemma \ref{lemma7}, we get three vertex disjoint paths $P(x^h,z^c),P(x^c,y^h)$ and $P(y^c,z^h)$ in $AQ_{n-1}^1.$ Similar to subcase 2.1 $Q_1,Q_2,\dots ,Q_{n-3},P_1,P_2,\dots ,P_{n-2}$ are the $2n-5$ vertex disjoint paths between $x$ and $y$ in $AQ_{n-2}^{00}.$
		Let $x_1,x_2,\dots ,x_{n-3}$ be the neighbours of x along $Q_1,Q_2,\dots ,Q_{n-3}$ respectively and $y_1,y_2,\dots ,y_{n-3}$be the neighbours of x along $Q_1,Q_2,\dots ,Q_{n-3}$ respectively. Also, $\{x_1^{h_{n-1}},x_2^{h_{n-1}},\dots ,x_{n-3}^{h_{n-1}}\}$ is the set of neighbours of $x^{h_{n-1}}$ in $AQ_{n-2}^{01}$ and $\{y_1^{h_{n-1}},y_2^{h_{n-1}},\dots ,y_{n-3}^{h_{n-1}}\}$ is the set of neighbours of $y^{h_{n-1}}$ in $AQ_{n-2}^{01}$. As $AQ_{n-2}^{00}$ is $2n-5$ connected, by Lemma \ref{lemma1}, there exist $2n-5$ vertex disjoint paths joining $z$ to every vertex in $\{x_1^{h_{n-1}},x_2^{h_{n-1}},\dots ,x_{n-3}^{h_{n-1}},y_1^{h_{n-1}},y_2^{h_{n-1}},\dots ,y_{n-4}^{h_{n-1}},y^{h_{n-1}},x^{h_{n-1}} \}$, say, $P(x_i^{h_{n-1}},z),i=1,2,\dots ,n-3$, $P(y_i^{h_{n-1}},z), i=1,2,\dots ,n-4 $, $P(y^{h_{n-1}},z)$ and $P(x^{h_{n-1}},z).$  we get $\psi_{n+\frac{n}{2}-6}$ paths with same procedure used in subcase $2.1$ and the remaining four paths as follows, (See Figure 8).\\
		$
		\psi_{n+\frac{n}{2}-5} = \{xx^{h_{n-1}}\}\cup P(x^{h_{n-1}}, z)  \cup P(y^{h_{n-1}}, z)\cup \{yy^{h_{n-1}}\} , \\ 
		\psi_{n+\frac{n}{2}-4} = P_{n-2} \cup \{x_{n-3}x_{n-3}^{h_{n-1}}\} \cup P( x_{n-3}^{h_{n-1}},z), \\ 
		\psi_{n+\frac{n}{2}-3} = P_{n-3}\cup \{xx^h\}\cup P(x^h, z^c) \cup \{z^cz\} \\ 
		\psi_{n+\frac{n}{2}-2} = \{xx^c\}\cup P(x^c, y^h) \cup \{y^hy,yy^c\} \cup P(y^c, z^h) \cup \{z^hz\}.
		$
		\\
		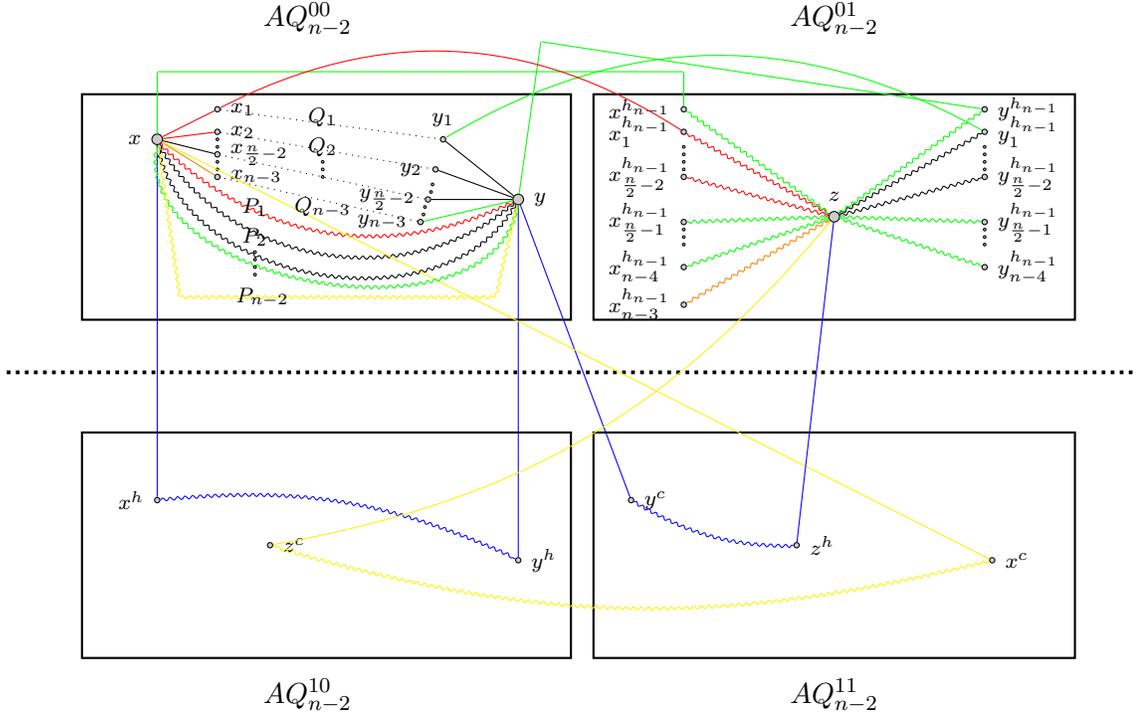
\begin{figure}[H]	
			\centering
			
			\begin{tikzpicture}
				% Draw the top-left rectangle
				\draw[thick] (-7, 2) rectangle (-0.5, 5);
				
				% Draw the top-right rectangle
				\draw[thick] (-0.2, 2) rectangle (6.2, 5);
				
				% Draw the bottom-left rectangle
				\draw[thick] (-7, -2.5) rectangle (-0.5, 0.5);
				
				% Draw the bottom-right rectangle
				\draw[thick] (-0.2, -2.5) rectangle (6.2, 0.5);
				
				% Labels for each rectangle
				\node at (-4, 6) {\textbf{$AQ_{n-2}^{00}$}};
				\node at (3, 6) {\textbf{$AQ_{n-2}^{01}$}};
				\node at (-4, -3) {\textbf{$AQ_{n-2}^{10}$}};
				\node at (3, -3) {\textbf{$AQ_{n-2}^{11}$}};
				
				% Dotted line
				\draw[dotted, line width=.5mm] (-8, 1.3) -- (7, 1.3);
				
				\node[circle, draw, fill=black!20, minimum size=4pt, inner sep=0pt, label={[font=\scriptsize]left:$x$}] (x) at (-6, 4.4) {};
				\node[circle, draw, fill=black!20, minimum size=2pt, inner  sep=0pt, label={[font=\scriptsize]left:$x^h$}] (xh) at (-6, -.4) {};
				\node[circle, draw, fill=black!20, minimum size=0pt, inner sep=0pt,] (b1) at (-6, 5.3) {};
				\draw[thin,green] (x) to (b1);
				\node[circle, draw, fill=black!20, minimum size=2pt, inner sep=0pt, label={[font=\scriptsize]right:$x_1$}] (x1) at (-5.2, 4.8) {};
				\node[circle, draw, fill=black!20, minimum size=2pt, inner sep=0pt, label={[font=\scriptsize]right:$x_2$}] (x2) at (-5.2, 4.5) {};
				\node[circle, draw, fill=black!20, minimum size=1pt, inner sep=0pt] (a) at (-5.2, 4.3) {};
				\node[circle, draw, fill=black!20, minimum size=1pt, inner sep=0pt] (a) at (-5.2, 4.4) {};
				\node[circle, draw, fill=black!20, minimum size=1pt, inner sep=0pt] (a) at (-5.2, 4.1) {};
				\node[circle, draw, fill=black!20, minimum size=1pt, inner sep=0pt] (a) at (-5.2, 4) {};
				\node[circle, draw, fill=black!20, minimum size=2pt, inner sep=0pt, label={[font=\scriptsize]right:$x_{n-3}$}] (xn-3) at (-5.2, 3.9) {};
				\node[circle, draw, fill=black!20, minimum size=2pt, inner sep=0pt, label={[font=\scriptsize]right:$x_{\frac{n}{2}-2}$}] (xn/2-2) at (-5.2, 4.2) {};
				\node[circle, draw, fill=black!20, minimum size=2pt, inner sep=0pt, label={[font=\scriptsize]above:$y_1$}] (y1) at (-2.2, 4.4) {};
				\node[circle, draw, fill=black!20, minimum size=2pt, inner sep=0pt, label={[font=\scriptsize]left:$y_2$}] (y2) at (-2.3, 4) {};
				\node[circle, draw, fill=black!20, minimum size=1pt, inner sep=0pt] (a) at (-2.34, 3.8) {};
				\node[circle, draw, fill=black!20, minimum size=1pt, inner sep=0pt] (a) at (-2.36, 3.7) {};
				\node[circle, draw, fill=black!20, minimum size=2pt, inner sep=0pt, label={[font=\scriptsize]left:$y_{n-3}$}] (yn-3) at (-2.5, 3.3) {};
				\node[circle, draw, fill=black!20, minimum size=1pt, inner sep=0pt] (a) at (-2.46, 3.4) {};
				\node[circle, draw, fill=black!20, minimum size=1pt, inner sep=0pt] (a) at (-2.44, 3.5) {};
				\node[circle, draw, fill=black!20, minimum size=2pt, inner sep=0pt, label={[font=\scriptsize]left:$y_{\frac{n}{2}-2}$}] (yn/2-2) at (-2.4, 3.6) {};
				\node[circle, draw, fill=black!20, minimum size=4pt, inner sep=0pt, label={[font=\scriptsize]right:$y$}] (y) at (-1.2, 3.6) {};
				\node[circle, draw, fill=black!20, minimum size=2pt, inner sep=0pt, label={[font=\scriptsize]right:$y^h$}] (yh) at (-1.2, -1.2) {};
				\node[circle, draw, fill=black!20, minimum size=2pt, inner sep=0pt, label={[font=\scriptsize]right:$x^c$}] (xc) at (5.1, -1.2) {};
				\node[circle, draw, fill=black!20, minimum size=2pt, inner sep=0pt, label={[font=\scriptsize]right:$y^c$}] (yc) at (.3, -.4) {};
				\node[circle, draw, fill=black!20, minimum size=2pt, inner sep=0pt, label={[font=\scriptsize]right:$z^h$}] (zh) at (2.5, -1) {};
				\node[circle, draw, fill=black!20, minimum size=2pt, inner sep=0pt, label={[font=\scriptsize]right:$z^c$}] (zc) at (-4.5, -1) {};
				\draw[thin,red] (x) to (x1);
				\draw[thin,red] (x) to (x2);
				\draw[decorate, decoration={snake, amplitude=0.02cm, segment length=0.1cm},,bend left=15,yellow] (xc) to (zc);
				\draw[decorate, decoration={snake, amplitude=0.02cm, segment length=0.1cm},bend left=15,blue] (xh) to (yh);
				\draw[decorate, decoration={snake, amplitude=0.02cm, segment length=0.1cm},bend right=15,blue] (yc) to (zh);	
				\draw[thin,orange] (x) to (xn-3);
				\draw[thin] (y) to (y1);
				\draw[thin] (y) to (y2);
				\draw[thin,green] (y) to (yn-3);
				\draw[dotted] (x1) to (y1);
				\draw[dotted] (x2) to (y2);
				\draw[dotted] (yn-3) to (xn-3);
				\draw[decorate, decoration={snake, amplitude=0.02cm, segment length=0.1cm},red, bend left=35] (y) to (x);
				\draw[decorate, decoration={snake, amplitude=0.02cm, segment length=0.1cm}, bend left=50] (y) to (x);
				\draw[thin] (x) to (xn/2-2);
				\draw[dotted] (yn/2-2) to (xn/2-2);
				\draw[thin] (yn/2-2) to (y);
				\draw[thin,blue] (x) to (xh);
				\draw[thin,blue] (y) to (yh);
				\draw[thin,blue] (y) to (yc);
				\draw[thin,yellow] (x) to (xc);
				
				\draw[decorate, decoration={snake, amplitude=0.02cm, segment length=0.1cm},green,bend left=90] (y) to (x);
				\draw[decorate, decoration={snake, amplitude=0.02cm, segment length=0.1cm}, bend left=70] (y) to (x);
				\node[circle, draw, fill=black!20, minimum size=0pt, inner sep=0pt] (a1) at (-5.7, 2.3) {};
				\node[circle, draw, fill=black!20, minimum size=0pt, inner sep=0pt] (a2) at (-1.5, 2.3) {};
				\node at (-4.7, 3.5) {\scriptsize{$P_1$}};
				\draw[decorate, decoration={snake, amplitude=0.02cm, segment length=0.1cm},yellow] (a1) to (x);
				\draw[decorate, decoration={snake, amplitude=0.02cm, segment length=0.1cm},yellow] (a2) to (y);
				\draw[decorate, decoration={snake, amplitude=0.02cm, segment length=0.1cm},yellow] (a1) to (a2);	
				\node at (-4.7, 3.1) {\scriptsize{$P_2$}};
				\node at (-4.6, 2.3) {\scriptsize{$P_{n-2}$}};
				\node at (-3.8, 4.7) {\scriptsize{$Q_1$}} ;
				\node at (-3.8, 4.3) {\scriptsize{$Q_2$}};
				\node at (-3.8, 3.5) {\scriptsize{$Q_{n-3}$}};
				\node[circle, draw, fill=black!20, minimum size=1pt, inner sep=0pt] (a) at (-3.8, 4.0) {};
				\node[circle, draw, fill=black!20, minimum size=1pt, inner sep=0pt] (a) at (-3.8, 4.1) {};
				\node[circle, draw, fill=black!20, minimum size=1pt, inner sep=0pt] (a) at (-3.8, 3.9) {};
				\node[circle, draw, fill=black!20, minimum size=1pt, inner sep=0pt] (a) at (-4.7, 2.9) {};
				\node[circle, draw, fill=black!20, minimum size=1pt, inner sep=0pt] (a) at (-4.7, 2.8) {};
				\node[circle, draw, fill=black!20, minimum size=1pt, inner sep=0pt] (a) at (-4.7, 2.7) {};
				\node[circle, draw, fill=black!20, minimum size=1pt, inner sep=0pt] (a) at (-4.7, 2.6) {};
				
				\node[circle, draw, fill=black!20, minimum size=2pt, inner sep=0pt, 
				label={[font=\scriptsize, align=left]left:$x_1^{h_{n-1}}$}] (x1hn-1) at (1, 4.5) {};
				\node[circle, draw, fill=black!20, minimum size=1pt, inner sep=0pt] (a) at (1, 4.2) {};
				\node[circle, draw, fill=black!20, minimum size=1pt, inner sep=0pt] (a) at (1, 4.1) {};
				\node[circle, draw, fill=black!20, minimum size=1pt, inner sep=0pt] (a) at (1, 4) {};
				\node[circle, draw, fill=black!20, minimum size=1pt, inner sep=0pt] (a) at (1, 4.3) {};
				\node[circle, draw, fill=black!20, minimum size=2pt, inner sep=0pt, 
				label={[font=\scriptsize, align=left]left:$x_{\frac{n}{2}-2}^{h_{n-1}}$}] (xn/2-2hn-1) at (1,3.9) {};
				\node[circle, draw, fill=black!20, minimum size=2pt, inner sep=0pt, 
				label={[font=\scriptsize, align=left]left:$x_{\frac{n}{2}-1}^{h_{n-1}}$}] (xn/2-1hn-1) at (1,3.3) {};
				\node[circle, draw, fill=black!20, minimum size=1pt, inner sep=0pt] (a) at (1, 3.2) {};
				\node[circle, draw, fill=black!20, minimum size=1pt, inner sep=0pt] (a) at (1, 3.1) {};
				\node[circle, draw, fill=black!20, minimum size=1pt, inner sep=0pt] (a) at (1, 3) {};
				\node[circle, draw, fill=black!20, minimum size=2pt, inner sep=0pt, 
				label={[font=\scriptsize, align=left]left:$x_{n-3}^{h_{n-1}}$}] (xn-3hn-1) at (1,2.2) {};
				\node[circle, draw, fill=black!20, minimum size=2pt, inner sep=0pt, 
				label={[font=\scriptsize, align=left]left:$x_{n-4}^{h_{n-1}}$}] (xn-4hn-1) at (1,2.7) {};
				
				\node[circle, draw, fill=black!20, minimum size=2pt, inner sep=0pt, 
				label={[font=\scriptsize, align=right]right:$y_1^{h_{n-1}}$}] (y1hn-1) at (5, 4.5) {};
				\node[circle, draw, fill=black!20, minimum size=1pt, inner sep=0pt] (a) at (5, 4.2) {};
				\node[circle, draw, fill=black!20, minimum size=1pt, inner sep=0pt] (a) at (5, 4.1) {};
				\node[circle, draw, fill=black!20, minimum size=1pt, inner sep=0pt] (a) at (5, 4) {};
				\node[circle, draw, fill=black!20, minimum size=1pt, inner sep=0pt] (a) at (5, 4.3) {};
				\node[circle, draw, fill=black!20, minimum size=2pt, inner sep=0pt, 
				label={[font=\scriptsize, align=right]right:$y_{\frac{n}{2}-2}^{h_{n-1}}$}] (yn/2-2hn-1) at (5,3.9) {};
				\node[circle, draw, fill=black!20, minimum size=2pt, inner sep=0pt, 
				label={[font=\scriptsize, align=right]right:$y_{\frac{n}{2}-1}^{h_{n-1}}$}] (yn/2-1hn-1) at (5,3.3) {};
				\node[circle, draw, fill=black!20, minimum size=1pt, inner sep=0pt] (a) at (5, 3.2) {};
				\node[circle, draw, fill=black!20, minimum size=1pt, inner sep=0pt] (a) at (5, 3.1) {};
				\node[circle, draw, fill=black!20, minimum size=1pt, inner sep=0pt] (a) at (5, 3) {};
				\node[circle, draw, fill=black!20, minimum size=2pt, inner sep=0pt, 
				label={[font=\scriptsize, align=right]right:$y_{n-4}^{h_{n-1}}$}] (yn-3hn-1) at (5,2.7) {};
				\node[circle, draw, fill=black!20, minimum size=4pt, inner sep=0pt, 
				label={[font=\scriptsize, ]above:$z$}] (z) at (3,3.37) {};
				\node[circle, draw, fill=black!20, minimum size=2pt, inner sep=0pt, 
				label={[font=\scriptsize, align=right]right:$y^{h_{n-1}}$}] (yhn-1) at (5,4.8) {};
				\node[circle, draw, fill=black!20, minimum size=2pt, inner sep=0pt, 
				label={[font=\scriptsize, align=left]left:$x^{h_{n-1}}$}] (xhn-1) at (1,4.8) {};
				\node[circle, draw, fill=black!20, minimum size=0pt, inner sep=0pt, ] (b2) at (1,5.3) {};
				\draw[thin,green] (xhn-1)to (b2);
				\draw[thin,green] (b1)to (b2);
				\draw[decorate, decoration={snake, amplitude=0.02cm, segment length=0.1cm}, thin,green] (z) -- (yn-3hn-1);
				\draw[decorate, decoration={snake, amplitude=0.02cm, segment length=0.1cm}, thin,green] (z) -- (yn/2-1hn-1);
				\draw[decorate, decoration={snake, amplitude=0.02cm, segment length=0.1cm}, thin,] (z) -- (y1hn-1);
				\draw[decorate, decoration={snake, amplitude=0.02cm, segment length=0.1cm}, thin, ] (z) -- (yn/2-2hn-1);
				\draw[decorate, decoration={snake, amplitude=0.02cm, segment length=0.1cm}, thin,orange] (z) -- (xn-3hn-1);
				\draw[decorate, decoration={snake, amplitude=0.02cm, segment length=0.1cm}, thin,green] (z) -- (xn/2-1hn-1);
				\draw[decorate, decoration={snake, amplitude=0.02cm, segment length=0.1cm}, thin, red] (z) -- (xn/2-2hn-1);
				\draw[decorate, decoration={snake, amplitude=0.02cm, segment length=0.1cm}, thin, red] (z) -- (x1hn-1);
				\draw[decorate, decoration={snake, amplitude=0.02cm, segment length=0.1cm}, thin, green] (z) -- (xhn-1);
				\draw[decorate, decoration={snake, amplitude=0.02cm, segment length=0.1cm}, thin, green] (z) -- (yhn-1);
				\draw[decorate, decoration={snake, amplitude=0.02cm, segment length=0.1cm}, thin,green] (z) -- (xn-4hn-1);
				\draw;	
				\node[circle, draw, fill=black!20, minimum size=0pt, inner sep=0pt, ] (b3) at (-0.9,5.7) {};
				\draw[thin,green, bend right=40] (y) -- (b3);
				\draw[thin,green, bend right=50] (yhn-1) -- (b3);
				\draw[thin,green,bend left=30] (y1)to (y1hn-1);
				\draw[thin,red,bend left=30] (x1)to(x1hn-1);
				\draw[thin,blue] (z) to (zh);
				\draw[thin,yellow,bend left=20] (z) to (zc);
			\end{tikzpicture}
			\caption{Case 2.2}
			\label{fig:graph-example}
		\end{figure}

		\textbf{Case 3 :} Let $x,y \in V(AQ_{n-1}^0)$ and $z\in V(AQ_{n-1}^1)$. As $AQ_{n-1}^0$ is $2n-3$ connected implies that there are $2n-3$ vertex disjoint paths joining $x$ and $y$ in $V(AQ_{n-1}^0)$. Without loss of generality, suppose $Q_1,Q_2,\dots ,Q_{n-2} $ are first $n-2$ paths among $2n-3$ paths we get. Let, $\{x_1,x_2,\dots ,x_{n-2}\}$ and $\{y_1,y_2,\dots,y_{n-2}\}$ are the sets of adjacent vertices of $x$ and $y$ respectively along $Q_1,Q_2,\dots ,Q_{n-2} $. 
		Consider, $X=\{x_1,x_2,\dots ,x_{n-2}, y_1,y_2,\dots,y_{n-3},x,y\}$ in $V(AQ_{n-1}^0)$ then, $X^h=\{x_1^h,x_2^h,\dots ,x_{n-2}^h, y_1^h,y_2^h,\dots,y_{n-3}^h,x^h,y^h\}$ in $V(AQ_{n-1}^1)$ is the set of $2n-3$ vertices which are adjacent to respective vertices of $X$.  As $AQ_{n-1}^1$ is $2n-3$ connected there exist $2n-3$ vertex disjoint paths joining $z$ to every vertex in $X^h$. Let $P(z,x_i^h),i=1,2,\dots ,n-2$; $P(z,y_i^h), i=1,2,\dots ,n-3; $ $P(z,x^h),$ and $P(z,y^h)$ be the paths joining $z$ to vertices in $X^h$. Let $P_1,P_2,\dots,P_{n-1}$ be the remaining $n-1$ paths between $x$ and $y$ in $AQ_{n-1}^{0}.$\\
		Then we get $2n-\frac{n}{2}-2$ vertex disjoint paths containing $x,y$ and $z$ in $AQ_n$ as follows, (See Figure 9).\\
		$
		\psi_1 = P_1 \cup \{xx_1\}\cup \{x_1x_1^{h_{n-1}}\} \cup P(x_1^{h_{n-1}}, z), \\ 
		\psi_2 = P_2 \cup \{xx_2\} \cup \{x_2x_2^{h_{n-1}}\} \cup P(x_2^{h_{n-1}}, z), \\ 
		\vdots \\ 
		\psi_{\frac{n}{2}-2} = P_{\frac{n}{2}-2} \cup \{xx_{\frac{n}{2}-2}\} \cup \{x_{\frac{n}{2}-2}x_{\frac{n}{2}-2}^{h_{n-1}}\} \cup P(x_{\frac{n}{2}-2}^{h_{n-1}}, z), \\ 
		\psi_{\frac{n}{2}-1} = P_{\frac{n}{2}-1}\cup \{yy_1\} \cup \{y_1y_1^{h_{n-1}}\} \cup P(y_1^{h_{n-1}}, z), \\ 
		\vdots \\ 
		\psi_{n-4} = P_{n-4}\cup \{yy_{\frac{n}{2}-2}\} \cup \{y_{\frac{n}{2}-2}y_{\frac{n}{2}-2}^{h_{n-1}}\} \cup P(y_{\frac{n}{2}-2}^{h_{n-1}}, z), \\ 
		\psi_{n-3} = \{x_{\frac{n}{2}-1}x_{\frac{n}{2}-1}^{h_{n-1}}\} \cup P(z, x_{\frac{n}{2}-1}^{h_{n-1}}) \cup \{y_{\frac{n}{2}-1}y_{\frac{n}{2}-1}^{h_{n-1}}\} \cup P(y_{\frac{n}{2}-1}^{h_{n-1}}, z), \\ 
		\vdots \\ 
		\psi_{n+\frac{n}{2}-5} = \{x_{n-3}x_{n-3}^{h_{n-1}}\} \cup P(z, x_{n-3}^{h_{n-1}}) \cup \{y_{n-3}y_{n-3}^{h_{n-1}}\} \cup P(y_{n-3}^{h_{n-1}}, z), \\ 
		\psi_{n+\frac{n}{2}-4} = \{yy^{h_{n-1}}\} \cup P(z, y^{h_{n-1}}) \cup \{zx\}, \\ 
		\psi_{n+\frac{n}{2}-3} = \{xx^h\} \cup \{zz^c\} \cup P_{n-3}, \\ 
		\psi_{n+\frac{n}{2}-2} = \{xx^c\} \cup \{zz^h\} \cup P_{n-2}.
		$
		\\
		\begin{figure} 
			\centering
			\begin{tikzpicture}
				% Draw the top-left rectangle
				\draw[thick] (-7, 1) rectangle (1, 6);
				
				% Draw the top-right rectangle
				\draw[thick] (-7, 0) rectangle (1, -5);

				% Labels for each rectangle
				\node at (-3, 6.5) {\textbf{$AQ_{n-1}^{0}$}};
				\node at (-3, -5.5) {\textbf{$AQ_{n-1}^{1}$}};

				\node[circle, draw, fill=black!20, minimum size=4pt, inner sep=0pt, label={[font=\scriptsize]left:$x$}] (x) at (-6, 3.3) {};
				\node[circle, draw, fill=black!20, minimum size=2pt, inner sep=0pt, label={[font=\scriptsize]right:$x_1$}] (x1) at (-5, 5.3) {};
				\node[circle, draw, fill=black!20, minimum size=2pt, inner sep=0pt, label={[font=\scriptsize]right:$x_2$}] (x2) at (-5, 4.9) {};
				\node[circle, draw, fill=black!20, minimum size=1pt, inner sep=0pt] (a) at (-5, 4.5) {};
				\node[circle, draw, fill=black!20, minimum size=1pt, inner sep=0pt] (a) at (-5, 4.3) {};
				\node[circle, draw, fill=black!20, minimum size=1pt, inner sep=0pt] (a) at (-5, 4.1) {};
				\node[circle, draw, fill=black!20, minimum size=2pt, inner sep=0pt, label={[font=\scriptsize]right:$x_{n-3}$}] (xn-1) at (-5, 3.5) {};
				\node[circle, draw, fill=black!20, minimum size=1pt, inner sep=0pt] (a) at (-5, 3.7) {};
				\node[circle, draw, fill=black!20, minimum size=2pt, inner sep=0pt, label={[font=\scriptsize]right:$x_{n-2}$}] (xn-3) at (-5, 3.3) {};
				\node[circle, draw, fill=black!20, minimum size=1pt, inner sep=0pt] (a) at (-5, 3.9) {};
				\node[circle, draw, fill=black!20, minimum size=2pt, inner sep=0pt, label={[font=\scriptsize]right:$x_{\frac{n}{2}-1}$}] (xn/2-2) at (-5, 4.1) {};
				\node[circle, draw, fill=black!20, minimum size=2pt, inner sep=0pt, label={[font=\scriptsize]left:$y_1$}] (y1) at (-1.2, 5.3) {};
				\node[circle, draw, fill=black!20, minimum size=2pt, inner sep=0pt, label={[font=\scriptsize]left:$y_2$}] (y2) at (-1.2, 4.9) {};
				\node[circle, draw, fill=black!20, minimum size=1pt, inner sep=0pt] (a) at (-1.2, 4.5) {};
				\node[circle, draw, fill=black!20, minimum size=1pt, inner sep=0pt] (a) at (-1.2, 4.3) {};
				\node[circle, draw, fill=black!20, minimum size=2pt, inner sep=0pt, label={[font=\scriptsize]left:$y_{n-2}$}] (yn-2) at (-1.2, 3.3) {};
				\node[circle, draw, fill=black!20, minimum size=1pt, inner sep=0pt] (a) at (-1.2, 3.9) {};
				\node[circle, draw, fill=black!20, minimum size=1pt, inner sep=0pt] (a) at (-1.2, 3.7) {};
				\node[circle, draw, fill=black!20, minimum size=1pt, inner sep=0pt] (a) at (-1.2, 3.5) {};
				\node[circle, draw, fill=black!20, minimum size=2pt, inner sep=0pt, label={[font=\scriptsize]left:$y_{\frac{n}{2}-1}$}] (yn/2-2) at (-1.2, 4.1) {};
				\node[circle, draw, fill=black!20, minimum size=4pt, inner sep=0pt, label={[font=\scriptsize]right:$y$}] (y) at (-0.2, 3.3) {};
				\node[circle, draw, fill=black!20, minimum size=2pt, inner sep=0pt, label={[font=\scriptsize]left:$y_{n-3}$}] (yn-3) at (-1.2, 3.5) {};

				\draw[thin, red] (x) to (x1);
				\draw[thin, red] (x) to (x2);
				\draw[thin,green] (x) to (xn-1);
				\draw[thin] (y) to (y1);
				\draw[thin] (y) to (y2);
				
				\draw[dotted] (x1) to (y1);
				\draw[dotted] (x2) to (y2);
				\draw[thick,blue] (x) to (xn-3);
				\draw[decorate, decoration={snake, amplitude=0.02cm, segment length=0.1cm}, red, bend left=15] (y) to (x);
				\draw[decorate, decoration={snake, amplitude=0.02cm, segment length=0.1cm}, red, bend left=30] (y) to (x);
				\draw[thin,red] (x) to (xn/2-2);
				\draw[dotted] (yn/2-2) to (xn/2-2);
				\draw[thin, bend right=8] (yn/2-2) to (y);
				\draw[decorate, decoration={snake, amplitude=0.02cm, segment length=0.1cm}, bend left=60] (y) to (x);
				\draw[decorate, decoration={snake, amplitude=0.02cm, segment length=0.1cm}, bend left=45] (y) to (x);
				\draw[decorate, decoration={snake, amplitude=0.02cm, segment length=0.1cm}, blue, bend left=75] (y) to (x);
				\node at (-3, 2.9) {\scriptsize $P_1$};
				\node at (-3, 1.5) {\scriptsize $P_{n-1}$};
				\node at (-3, 1.8) {\scriptsize $P_{n-2}$};
				\node[circle, draw, fill=black!20, minimum size=1pt, inner sep=0pt] (a) at (-3, 2.7) {};
				\node[circle, draw, fill=black!20, minimum size=1pt, inner sep=0pt] (a) at (-3, 2.6) {};
				\node[circle, draw, fill=black!20, minimum size=1pt, inner sep=0pt] (a) at (-3, 2.3) {};
				\node[circle, draw, fill=black!20, minimum size=1pt, inner sep=0pt] (a) at (-3, 2.2) {};
				
				\node at (-3, 5.3) {\scriptsize{$Q_1$}} ;
				\node at (-3, 4.9) {\scriptsize{$Q_2$}};
				\node at (-3, 3.3) {\scriptsize{$Q_{n-2}$}};
				
				\draw[thin] (y) to (yn-2);
				\draw[thin,green] (y) to (yn-3);
				\draw[dotted] (xn-1) to (yn-3);
				\draw[dotted] (xn-3) to (yn-2);
				\node[circle, draw, fill=black!20, minimum size=2pt, inner sep=0pt, 
				label={[font=\scriptsize, align=left]left:$x_1^h$}] (x1hn-1) at (-5, -1.5) {};
				\node[circle, draw, fill=black!20, minimum size=1pt, inner sep=0pt] (a) at (-5, -1.7) {};
				\node[circle, draw, fill=black!20, minimum size=1pt, inner sep=0pt] (a) at (-5, -1.9) {};
				\node[circle, draw, fill=black!20, minimum size=1pt, inner sep=0pt] (a) at (-5, -3.2) {};
				
				\node[circle, draw, fill=black!20, minimum size=2pt, inner sep=0pt, 
				label={[font=\scriptsize, align=left]left:$x_{\frac{n}{2}-2}^h$}] (xn/2-2hn-1) at (-5, -2.1) {};
				\node[circle, draw, fill=black!20, minimum size=2pt, inner sep=0pt, 
				label={[font=\scriptsize, align=left]left:$x_{\frac{n}{2}-1}^h$}] (xn/2-1hn-1) at (-5, -2.6) {};
				
				\node[circle, draw, fill=black!20, minimum size=1pt, inner sep=0pt] (a) at (-5, -3.4) {};
				
				\node[circle, draw, fill=black!20, minimum size=2pt, inner sep=0pt, 
				label={[font=\scriptsize, align=left]left:$x_{n-3}^h$}] (xn-3hn-1) at (-5, -3.6) {};
				\node[circle, draw, fill=black!20, minimum size=2pt, inner sep=0pt, 
				label={[font=\scriptsize, align=left]left:$x_{n-2}^h$}] (xn-3hn-2) at (-5, -4) {};
				
				\node[circle, draw, fill=black!20, minimum size=2pt, inner sep=0pt, 
				label={[font=\scriptsize, align=right]right:$y_1^h$}] (y1hn-1) at (-1.2, -1.5) {};
				\node[circle, draw, fill=black!20, minimum size=1pt, inner sep=0pt] (a) at (-1.2, -1.7) {};
				\node[circle, draw, fill=black!20, minimum size=1pt, inner sep=0pt] (a) at (-1.2, -1.9) {};
				\node[circle, draw, fill=black!20, minimum size=1pt, inner sep=0pt] (a) at (-1.2, -3.2) {};
				\node[circle, draw, fill=black!20, minimum size=1pt, inner sep=0pt] (a) at (-1.2, -3.4) {};
				\node[circle, draw, fill=black!20, minimum size=2pt, inner sep=0pt, 
				label={[font=\scriptsize, align=right]right:$y_{\frac{n}{2}-2}^h$}] (yn/2-2hn-1) at (-1.20, -2.1) {};
				\node[circle, draw, fill=black!20, minimum size=2pt, inner sep=0pt, 
				label={[font=\scriptsize, align=right]right:$y_{\frac{n}{2}-1}^h$}] (yn/2-1hn-1) at (-1.20, -2.6) {};
				\node[circle, draw, fill=black!20, minimum size=2pt, inner sep=0pt, 
				label={[font=\scriptsize, align=right]right:$y_{n-3}^h$}] (yn-3hn-1) at (-1.20, -3.6) {};
				\node[circle, draw, fill=black!20, minimum size=4pt, inner sep=0pt, 
				label=above:$z$] (z) at (-3, -3) {};
				\node[circle, draw, fill=black!20, minimum size=4pt, inner sep=0pt, 
				label={[font=\scriptsize, align=right]right:$y^h$}] (yhn-1) at (0, -3) {};
				\node[circle, draw, fill=black!20, minimum size=4pt, inner sep=0pt, 
				label={[font=\scriptsize, align=left]left:$x^h$}] (xhn-1) at (-6, -3) {};

				\draw[decorate, decoration={snake, amplitude=0.02cm, segment length=0.1cm}, thin, green] (z) -- (yn-3hn-1);
				\draw[decorate, decoration={snake, amplitude=0.02cm, segment length=0.1cm}, thin] (z) -- (yn/2-1hn-1);
				\draw[decorate, decoration={snake, amplitude=0.02cm, segment length=0.1cm}, thin] (z) -- (y1hn-1);
				\draw[decorate, decoration={snake, amplitude=0.02cm, segment length=0.1cm}, thin] (z) -- (yn/2-2hn-1);
				\draw[decorate, decoration={snake, amplitude=0.02cm, segment length=0.1cm}, thin, green] (z) -- (xn-3hn-1);
				\draw[decorate, decoration={snake, amplitude=0.02cm, segment length=0.1cm}, thin, red] (z) -- (xn/2-1hn-1);
				\draw[decorate, decoration={snake, amplitude=0.02cm, segment length=0.1cm}, thin, red] (z) -- (xn/2-2hn-1);
				\draw[decorate, decoration={snake, amplitude=0.02cm, segment length=0.1cm}, thin, red] (z) -- (x1hn-1);
				\draw[decorate, decoration={snake, amplitude=0.02cm, segment length=0.1cm}, thin, green] (z) -- (xhn-1);
				\draw[decorate, decoration={snake, amplitude=0.02cm, segment length=0.1cm}, thin, green] (z) -- (yhn-1);
				\draw[decorate, decoration={snake, amplitude=0.02cm, segment length=0.1cm}, thin, blue] (z) -- (xn-3hn-2);

				\draw[thin, red, bend right=60] (x1) to (x1hn-1);
				%	\draw[thin, bend left=50,red] (xn/2-2) to (xn/2-1hn-1);
				\draw[thin, bend left=60] (y1) to (y1hn-1);
				%%	\draw[thin, bend left=40] (yn/2-2) to (yn/2-1hn-1);
				\draw[thin, bend left=40, blue] (xn-3) to (xn-3hn-2);
				%	\draw[thick,blue] (y) to (yn-3);
				\draw[thin,green, bend right=30] (yn-3) to (yn-3hn-1);
				
				\draw[thin,green, bend right=30] (xn-3hn-1) to (xn-1);
				\draw[dotted, line width=.5mm] (-7.5, .5) -- (1.5, .5);
			\end{tikzpicture}
			
			\caption{Case 3}
			\label{fig:graph-example}
		\end{figure}
	\end{proof}
	Let us prove the result for odd integer $n.$
	\begin{corollary}
		$\pi_3(AQ_n)=\frac{3(n-1)}{2}-1$, when $n$ is odd. 
	\end{corollary}
	\begin{proof}
		Since $n$ is odd, by Lemma \ref{lemma4}, $\pi_3(AQ_n)\leq \frac{3(n-1)}{2}-1$. Now we have to prove that $\pi_3(AQ_n)= \frac{3(n-1)}{2}-1$ by constructing $\frac{3(n-1)}{2}-1$ vertex disjoint paths containing any three distinct vertices in $V(AQ_n).$ Let, \{x,y,z\} be any three distinct vertices of $AQ_n.$\\
		\textbf{Case: 1} Suppose, $\{x,y,z\}\subseteq V(AQ^i_{n-1})$ for some $0\leq i\leq 1.$ Without loss of generality assume that $\{x,y,z\} \subseteq AQ_{n-1}^0$, then $x^h,y^h,z^h,x^c\in AQ_{n-1}^1.$
		As $n$ is an odd integer, $n-1$ is even. Thus, by Theorem $1$, there are $\frac{3(n-1)}{2}-2$ vertex disjoint paths containing $x,y$ and $z$ in $AQ_{n-1}.$ Thus we have to find only one additional path containing $\{x,y,z\}$ which is vertex disjoint from earlier obtained $\frac{3(n-1)}{2}-2$ paths. Now, by Lemma \ref{lemma7} there exist two disjoint paths $P(x^h,y^h)$ and $P(x^c,z^h)$ in $AQ_{n-1}^1.$ \\
		We get the required remaining path containing $x,y$ and $z$ in $AQ_n$ as follows, (See Figure 10).\\
		$\psi_{\frac{3(n-1)}{2}-1}=
		\{yy^h\} \cup P(x^h, y^h)\cup \{xx^h\} \cup \{xx^c\} \cup P(x^c, z^h)\cup \{zz^h\}. $\\
		
		\textbf{Case: 2} Suppose $\{x,y\}\subseteq V(AQ_{n-1}^0)$ and $z\in V(AQ_{n-1}^1) $\\
		Since $AQ_{n-1}$ is $2n-3$ connected, there exist $2n-3$ paths joining $x$ and $y$ in $AQ_{n-1}^0$. Let $Q_1, Q_2, \dots, Q_{n-2}$ be the paths between $x$ and $y$. Denote the adjacent vertices of $x$ and $y$ along these paths by $x_1, x_2, \dots, x_{n-2}$ and $y_1, y_2, \dots, y_{n-2}$, respectively. Let, $P_1,P_2,\dots P_{n-1}$ be the remaining paths between $x$ and $y$ in $AQ_{n-1}^0$

		\begin{figure}[h]
			\centering
			\begin{tikzpicture}
				\draw[thick] (-7, 1) rectangle (1, 6);
				
				% Draw the top-right rectangle
				\draw[thick] (-7, 0) rectangle (1, -5);
				
				\draw[dotted, line width=.5mm] (-7.5, .5) -- (1.5, .5);
				% Labels for each rectangle
				\node at (-3, 6.5) {\textbf{$AQ_{n-1}^{0}$}};
				\node at (-3, -5.5) {\textbf{$AQ_{n-1}^{1}$}};
				% Nodes and labels
				\node[circle, draw, fill=black!20, minimum size=4pt, inner sep=0pt, label={[font=\scriptsize]left:$x$}] (x) at (-6, 5) {};
				\node[circle, draw, fill=black!20, minimum size=4pt, inner sep=0pt, label={[font=\scriptsize]left:$z$}] (z) at (-2, 1.5) {};
				\node[circle, draw, fill=black!20, minimum size=4pt, inner sep=0pt, label={[font=\scriptsize]left:$y$}] (y) at (-4, 3.2) {};
				\node[circle, draw, fill=black!20, minimum size=4pt, inner sep=0pt, label={[font=\scriptsize]right:$x^h$}] (xh) at (-6, -1) {};
				\node[circle, draw, fill=black!20, minimum size=4pt, inner sep=0pt, label={[font=\scriptsize]right:$x^c$}] (xc) at (-1, -4) {};
				\node[circle, draw, fill=black!20, minimum size=4pt, inner sep=0pt, label={[font=\scriptsize]left:$z^h$}] (zh) at (-2, -4) {};
				\node[circle, draw, fill=black!20, minimum size=4pt, inner sep=0pt, label={[font=\scriptsize]right:$y^h$}] (yh) at (-4, -2) {};
				
				% Edges
				\draw[decorate, decoration={snake, amplitude=0.02cm, segment length=0.1cm},red] (xh) -- (yh);
				\draw[decorate, decoration={snake, amplitude=0.02cm, segment length=0.1cm},red] (xc) -- (zh); 
				\draw[thin,red] (xh) -- (x);
				\draw[thin,red] (yh) -- (y);
				\draw[thin,red] (zh) -- (z); 
				\draw[thin,red] (x) -- (xc); 
			\end{tikzpicture}
			\caption{Case 1}
		\end{figure}
		
		Now consider the set $X = \{x_1, x_2, \dots, x_{n-3}, y_1, y_2, \dots, y_{n-3}, x, y\}\subseteq V(AQ^0_{n-1})$, then $X^h= \{x_1^h, x_2^h, \dots, x_{n-3}^h, y_1^h, y_2^h, \dots, y_{n-3}^h, x^h, y^h\}$ is the set of adjacent vertices of $X$ in $AQ_{n-1}^1$. Since $AQ_{n-1}^1$ is $2n-3$ connected there exist $2n-3$ vertex-disjoint paths in $AQ_{n-1}^0$ joining $z$ to each vertex in $X^h$, Say, $
		P(x_1^h, z), P(x_2^h, z), \dots, P(x_{n-3}^h, z), P(y_1^h, z), P(y_2^h, z),$  $\dots, P(y_{n-3}^h, z), P(x^h, z), P(y^h, z)$.
		
		We obtain a total of $\frac{3(n-1)}{2} - 1$ paths in $AQ_n$ containing $x$, $y$, and $z$ as follows, (See Figure 11).\\
		$
		\psi_i =\{xx_i\}\cup \{x_ix_i^h\} \cup P(x_i^h, z) \cup P_i,1\leq i\leq \frac{n-1}{2} \\ 
		\psi_{\frac{n+1}{2}} = \{yy_1\}\cup \{y_1y_1^h\} \cup P(y_1^h, z) \cup P_{\frac{n+1}{2}}, \\ 
		\vdots \\ 
		\psi_{n-1} = \{yy_{\frac{n-1}{2}}\} \cup \{y_{\frac{n-1}{2}}y_{\frac{n-1}{2}}^h\} \cup P(y_{\frac{n-1}{2}}^h, z) \cup P_{n-1}, \\ 
		\psi_n = \{xx_{\frac{n+1}{2}}\} \cup \{x_{\frac{n+1}{2}}x_{\frac{n+1}{2}}^h\} \cup P(x_{\frac{n+1}{2}}^h, z)\cup P(y_{\frac{n+1}{2}}^h,z)
		\cup \{yy_{\frac{n+1}{2}}\} \cup \{y_{\frac{n+1}{2}}y_{\frac{n+1}{2}}^h\}  , \\ 
		\vdots \\ 
		\psi_{\frac{3(n-1)}{2}-2} = \{xx_{n-3}\} \cup \{x_{n-3}x_{n-3}^h\} \cup P(x_{n-3}^h, z)\cup P(y_{n-3}^h,z)\cup \{y_{n-3}y_{n-3}^h\} \cup \{yy_{n-3}\}   , \\ 
		\psi_{\frac{3(n-1)}{2}-1} = \{xx^h\}\cup P(x^h, z) \cup P(y^h, z)\cup \{yy^h\}  .
		$

		\begin{figure} [H]
			\centering
			\begin{tikzpicture}
				% Draw the top-left rectangle
				\draw[thick] (-7, 1) rectangle (1, 6);
				
				% Draw the top-right rectangle
				\draw[thick] (-7, 0) rectangle (1, -5);
				\draw[dotted, line width=.5mm] (-7.5, .5) -- (1.5, .5);
				
				% Labels for each rectangle
				\node at (-3, 6.5) {\textbf{$AQ_{n-1}^{0}$}};
				\node at (-3, -5.5) {\textbf{$AQ_{n-1}^{1}$}};

				\node[circle, draw, fill=black!20, minimum size=4pt, inner sep=0pt, label={[font=\scriptsize]left:$x$}] (x) at (-6, 3.3) {};
				\node[circle, draw, fill=black!20, minimum size=2pt, inner sep=0pt, label={[font=\scriptsize]right:$x_1$}] (x1) at (-5, 5.3) {};
				\node[circle, draw, fill=black!20, minimum size=2pt, inner sep=0pt, label={[font=\scriptsize]right:$x_2$}] (x2) at (-5, 4.9) {};
				\node[circle, draw, fill=black!20, minimum size=1pt, inner sep=0pt] (a) at (-5, 4.5) {};
				\node[circle, draw, fill=black!20, minimum size=1pt, inner sep=0pt] (a) at (-5, 4.3) {};
				\node[circle, draw, fill=black!20, minimum size=1pt, inner sep=0pt] (a) at (-5, 4.1) {};
				\node[circle, draw, fill=black!20, minimum size=2pt, inner sep=0pt, label={[font=\scriptsize]right:$x_{n-3}$}] (xn-1) at (-5, 3.5) {};
				\node[circle, draw, fill=black!20, minimum size=1pt, inner sep=0pt] (a) at (-5, 3.7) {};
				\node[circle, draw, fill=black!20, minimum size=2pt, inner sep=0pt, label={[font=\scriptsize]right:$x_{n-2}$}] (xn-3) at (-5, 3.3) {};
				\node[circle, draw, fill=black!20, minimum size=1pt, inner sep=0pt] (a) at (-5, 3.9) {};
				\node[circle, draw, fill=black!20, minimum size=2pt, inner sep=0pt, label={[font=\scriptsize]right:$x_{\frac{n-1}{2}}$}] (xn/2-2) at (-5, 4.1) {};
				\node[circle, draw, fill=black!20, minimum size=2pt, inner sep=0pt, label={[font=\scriptsize]left:$y_1$}] (y1) at (-1.2, 5.3) {};
				\node[circle, draw, fill=black!20, minimum size=2pt, inner sep=0pt, label={[font=\scriptsize]left:$y_2$}] (y2) at (-1.2, 4.9) {};
				\node[circle, draw, fill=black!20, minimum size=1pt, inner sep=0pt] (a) at (-1.2, 4.5) {};
				\node[circle, draw, fill=black!20, minimum size=1pt, inner sep=0pt] (a) at (-1.2, 4.3) {};
				\node[circle, draw, fill=black!20, minimum size=2pt, inner sep=0pt, label={[font=\scriptsize]left:$y_{n-2}$}] (yn-2) at (-1.2, 3.3) {};
				\node[circle, draw, fill=black!20, minimum size=1pt, inner sep=0pt] (a) at (-1.2, 3.9) {};
				\node[circle, draw, fill=black!20, minimum size=1pt, inner sep=0pt] (a) at (-1.2, 3.7) {};
				\node[circle, draw, fill=black!20, minimum size=1pt, inner sep=0pt] (a) at (-1.2, 3.5) {};
				\node[circle, draw, fill=black!20, minimum size=2pt, inner sep=0pt, label={[font=\scriptsize]left:$y_{\frac{n-1}{2}}$}] (yn/2-2) at (-1.2, 4.1) {};
				\node[circle, draw, fill=black!20, minimum size=4pt, inner sep=0pt, label={[font=\scriptsize]right:$y$}] (y) at (-0.2, 3.3) {};
				\node[circle, draw, fill=black!20, minimum size=2pt, inner sep=0pt, label={[font=\scriptsize]left:$y_{n-3}$}] (yn-3) at (-1.2, 3.5) {};

				\draw[thin, red] (x) to (x1);
				\draw[thin, red] (x) to (x2);
				\draw[thin,green] (x) to (xn-1);
				\draw[thin] (y) to (y1);
				\draw[thin] (y) to (y2);
				
				\draw[dotted] (x1) to (y1);
				\draw[dotted] (x2) to (y2);
				\draw[thin] (x) to (xn-3);
				\draw[decorate, decoration={snake, amplitude=0.02cm, segment length=0.1cm}, red, bend left=15] (y) to (x);
				\draw[decorate, decoration={snake, amplitude=0.02cm, segment length=0.1cm}, red, bend left=30] (y) to (x);
				\draw[thin,red] (x) to (xn/2-2);
				\draw[dotted] (yn/2-2) to (xn/2-2);
				\draw[thin, bend right=8] (yn/2-2) to (y);
				\draw[decorate, decoration={snake, amplitude=0.02cm, segment length=0.1cm}, bend left=60] (y) to (x);
				\draw[decorate, decoration={snake, amplitude=0.02cm, segment length=0.1cm}, bend left=45] (y) to (x);
				\draw[decorate, decoration={snake, amplitude=0.02cm, segment length=0.1cm}, bend left=75] (y) to (x);
				\node at (-3, 2.9) {\scriptsize $P_1$};
				\node at (-3, 1.5) {\scriptsize $P_{n-1}$};
				\node at (-3, 1.8) {\scriptsize $P_{n-2}$};
				\node[circle, draw, fill=black!20, minimum size=1pt, inner sep=0pt] (a) at (-3, 2.7) {};
				\node[circle, draw, fill=black!20, minimum size=1pt, inner sep=0pt] (a) at (-3, 2.6) {};
				\node[circle, draw, fill=black!20, minimum size=1pt, inner sep=0pt] (a) at (-3, 2.3) {};
				\node[circle, draw, fill=black!20, minimum size=1pt, inner sep=0pt] (a) at (-3, 2.2) {};
				
				\node at (-3, 5.3) {\scriptsize{$Q_1$}} ;
				\node at (-3, 4.9) {\scriptsize{$Q_2$}};
				\node at (-3, 3.3) {\scriptsize{$Q_{n-2}$}};
				
				\draw[thin] (y) to (yn-2);
				\draw[thin,green] (y) to (yn-3);
				\draw[dotted] (xn-1) to (yn-3);
				\draw[dotted] (xn-3) to (yn-2);
				\node[circle, draw, fill=black!20, minimum size=2pt, inner sep=0pt, 
				label={[font=\scriptsize, align=left]left:$x_1^h$}] (x1hn-1) at (-5, -1.5) {};
				\node[circle, draw, fill=black!20, minimum size=1pt, inner sep=0pt] (a) at (-5, -1.7) {};
				\node[circle, draw, fill=black!20, minimum size=1pt, inner sep=0pt] (a) at (-5, -1.9) {};
				\node[circle, draw, fill=black!20, minimum size=1pt, inner sep=0pt] (a) at (-5, -3.2) {};
				
				\node[circle, draw, fill=black!20, minimum size=2pt, inner sep=0pt, 
				label={[font=\scriptsize, align=left]left:$x_{\frac{n-1}{2}}^h$}] (xn/2-2hn-1) at (-5, -2.1) {};
				\node[circle, draw, fill=black!20, minimum size=2pt, inner sep=0pt, 
				label={[font=\scriptsize, align=left]left:$x_{\frac{n+1}{2}}^h$}] (xn/2-1hn-1) at (-5, -2.6) {};
				
				\node[circle, draw, fill=black!20, minimum size=1pt, inner sep=0pt] (a) at (-5, -3.4) {};
				
				\node[circle, draw, fill=black!20, minimum size=2pt, inner sep=0pt, 
				label={[font=\scriptsize, align=left]left:$x_{n-3}^h$}] (xn-3hn-1) at (-5, -3.6) {};

				\node[circle, draw, fill=black!20, minimum size=2pt, inner sep=0pt, 
				label={[font=\scriptsize, align=right]right:$y_1^h$}] (y1hn-1) at (-1.2, -1.5) {};
				\node[circle, draw, fill=black!20, minimum size=1pt, inner sep=0pt] (a) at (-1.2, -1.7) {};
				\node[circle, draw, fill=black!20, minimum size=1pt, inner sep=0pt] (a) at (-1.2, -1.9) {};
				\node[circle, draw, fill=black!20, minimum size=1pt, inner sep=0pt] (a) at (-1.2, -3.2) {};
				\node[circle, draw, fill=black!20, minimum size=1pt, inner sep=0pt] (a) at (-1.2, -3.4) {};
				\node[circle, draw, fill=black!20, minimum size=2pt, inner sep=0pt, 
				label={[font=\scriptsize, align=right]right:$y_{\frac{n-1}{2}}^h$}] (yn/2-2hn-1) at (-1.20, -2.1) {};
				\node[circle, draw, fill=black!20, minimum size=2pt, inner sep=0pt, 
				label={[font=\scriptsize, align=right]right:$y_{\frac{n+1}{2}}^h$}] (yn/2-1hn-1) at (-1.20, -2.6) {};
				\node[circle, draw, fill=black!20, minimum size=2pt, inner sep=0pt, 
				label={[font=\scriptsize, align=right]right:$y_{n-3}^h$}] (yn-3hn-1) at (-1.20, -3.6) {};
				\node[circle, draw, fill=black!20, minimum size=4pt, inner sep=0pt, 
				label=above:$z$] (z) at (-3, -3) {};
				\node[circle, draw, fill=black!20, minimum size=4pt, inner sep=0pt, 
				label={[font=\scriptsize, align=right]right:$y^h$}] (yhn-1) at (0, -3) {};
				\node[circle, draw, fill=black!20, minimum size=4pt, inner sep=0pt, 
				label={[font=\scriptsize, align=left]left:$x^h$}] (xhn-1) at (-6, -3) {};

				\draw[decorate, decoration={snake, amplitude=0.02cm, segment length=0.1cm}, thin, green] (z) -- (yn-3hn-1);
				\draw[decorate, decoration={snake, amplitude=0.02cm, segment length=0.1cm}, thin, green] (z) -- (yn/2-1hn-1);
				\draw[decorate, decoration={snake, amplitude=0.02cm, segment length=0.1cm}, thin] (z) -- (y1hn-1);
				\draw[decorate, decoration={snake, amplitude=0.02cm, segment length=0.1cm}, thin] (z) -- (yn/2-2hn-1);
				\draw[decorate, decoration={snake, amplitude=0.02cm, segment length=0.1cm}, thin, green] (z) -- (xn-3hn-1);
				\draw[decorate, decoration={snake, amplitude=0.02cm, segment length=0.1cm}, thin, green] (z) -- (xn/2-1hn-1);
				\draw[decorate, decoration={snake, amplitude=0.02cm, segment length=0.1cm}, thin, red] (z) -- (xn/2-2hn-1);
				\draw[decorate, decoration={snake, amplitude=0.02cm, segment length=0.1cm}, thin, red] (z) -- (x1hn-1);
				\draw[decorate, decoration={snake, amplitude=0.02cm, segment length=0.1cm}, thin, green] (z) -- (xhn-1);
				\draw[decorate, decoration={snake, amplitude=0.02cm, segment length=0.1cm}, thin, green] (z) -- (yhn-1);

				\draw[thin, red, bend right=60] (x1) to (x1hn-1);
				%	\draw[thin, bend left=50,red] (xn/2-2) to (xn/2-1hn-1);
				\draw[thin, bend left=60] (y1) to (y1hn-1);
				%%	\draw[thin, bend left=40] (yn/2-2) to (yn/2-1hn-1);
				
				%	\draw[thick,blue] (y) to (yn-3);
				\draw[thin,green, bend right=30] (yn-3) to (yn-3hn-1);
				
				\draw[thin,green, bend right=30] (xn-3hn-1) to (xn-1);
			\end{tikzpicture}
			
			\caption{Case 2}
			\label{fig:graph-example}
		\end{figure}
	\end{proof}

	\section*{Acknowledgment}
	
	The Third (corresponding) author gratefully acknowledges the Department of Science and Technology, New Delhi, India, for awarding the Women Scientist Scheme (DST/WOS-A/PM-14/2021(G)) for research in Basic/Applied Sciences.


\begin{thebibliography}{99}
		\bibitem{chartrand1984generalized} G.~Chartrand, F.~Kapoor, L.~Lesniak, and D.~R.~Lick, ``Generalized connectivity in graphs,'' \textit{Bombay Math. Colloq. Bull.}, $2$ ($1984$), $1$--$6$.
		
		\bibitem{choudum} S. A. Choudum, V. Sunitha, Augmented cubes, \textit{Networks: An International Journal}, $40$ ($2002$) $71$--$84$.
		
		\bibitem{hager} M. Hager, Path-connectivity in graphs, \textit{Discrete Mathematics}, $59$ ($1986$) $53$--$59$.
		
		\bibitem{li2021tree} S. Li, Z. Qin, J. Tu, J. Yue, On tree-connectivity and path-connectivity of graphs, \textit{Graphs and Combinatorics}, $37$ ($2021$) $2521$--$2533$.
		
		\bibitem{mao} Y. Mao, Path-connectivity of lexicographic product graphs, \textit{International Journal of Computer Mathematics}, $93$ ($2016$) $27$--$39$.
		
		\bibitem{C1} M. Ma, G. Liu, J.-M. Xu, The super connectivity of augmented cubes, \textit{Information Processing Letters}, $106$ ($2008$) $59$--$63$.
		
		\bibitem{wang} Wang, J., Cheng, D. The $3$-path-connectivity of pancake graphs, \textit{The Journal of Supercomputing}, $81$ ($2025$) $1$–$30$.
		
		\bibitem{west2001} D. B. West, \textit{Introduction to Graph Theory}, $2$nd edn., Pearson Education, Delhi, $2001$.
		
		\bibitem{zhu} W.-H. Zhu, R.-X. Hao, L. Li, The $3$-path-connectivity of the hypercubes, \textit{Discrete Applied Mathematics}, $322$ ($2022$) $203$--$209$.
		
		\bibitem{C2} W.-H. Zhu, R.-X. Hao, Y.-Q. Feng, J. Lee, The $3$-path-connectivity of the $k$-ary $n$-cube, \textit{Applied Mathematics and Computation}, $436$ ($2023$) $127499$.
		
		\bibitem{generalizedhypercube} S. Lin, Q. Zhang, The generalized $4$-connectivity of hypercubes, \textit{Discrete Applied Mathematics}, $220$ ($2017$) $60$--$67$.
		
		\bibitem{generalizedfolded} H. Liu, D. Cheng, The generalized $4$-connectivity of folded hypercube, \textit{International Journal of Computer Mathematics: Computer Systems Theory}, $7$ ($2022$) $235$--$245$.
		
		\bibitem{mane2024pendant} S. A. Mane, S. A. Kandekar, Pendant $3$-tree-connectivity of augmented cubes, \textit{The Journal of Supercomputing}, $80$ ($2024$) $19395$--$19413$.
		
		\bibitem{wang20253} Y. Wang, D. Cheng, The $3$-path connectivity of the folded hypercube, \textit{The Journal of Supercomputing}, $81$ ($2025$) $1$--$25$.
		
		\bibitem{cheng2023generalized} D. Cheng, The generalized $4$-connectivity of locally twisted cubes, \textit{Journal of Applied Mathematics and Computing}, $69$ ($2023$) $3095$--$3111$.
		
		%\bibitem{kandekar2020} S. A. Kandekar, On fault tolerance in graphical structures and related aspects, PhD Thesis, Savitribai Phule Pune University, $2020$.		
	\end{thebibliography}
\end{document}